\documentclass[11pt]{article}
\usepackage{lineno}
\usepackage[margin=1in]{geometry}
\usepackage{amsmath,amsfonts,amsthm,amssymb}

\usepackage{cite}
\usepackage{subcaption}
\usepackage{graphicx}
\graphicspath{{figures/}}
\usepackage{color}
\usepackage{url}
\usepackage{soul}
\usepackage{float}
\usepackage{enumitem}
\usepackage{hhline}
\usepackage{tabularx}
\usepackage{array}
\usepackage{calc}
\usepackage{tikz}

\usepackage{todonotes}

\usepackage{verbatim}

\newcommand{\bbR}{\mathbb{R}}

\newcommand{\Rn}{\mathbb{R}^n}

\newcommand{\collar}{\Omega_I}

\newcommand{\Ldel}{\mathcal{L}^{\delta,\beta}}

\newcommand{\cdel}{c^{\delta,\beta}}
\newcommand{\mdel}{m^{\delta,\beta}}

\newtheorem{conjecture}{Conjecture}

\theoremstyle{definition}
\newtheorem{ex}{Example}

\theoremstyle{remark}
\newtheorem{remark}{Remark}

\title{Fourier Spectral Method for Nonlocal Equations on Bounded Domains
  \thanks{This project is based upon work supported by the National Science Foundation under Grant No. 2108588.}
}

\author{%
  Ilyas Mustapha\thanks{ilyas1@ksu.edu, Department of Mathematics, Kansas State University, Manhattan, KS.} \quad
  Bacim Alali\thanks{bacimalali@math.ksu.edu, Department of Mathematics, Kansas State University, Manhattan, KS.} \quad
  Nathan Albin\thanks{albin@ksu.edu, Department of Mathematics, Kansas State University, Manhattan, KS.}
}
\date{}
\begin{document}        
\maketitle

\begin{abstract}
This work introduces efficient and accurate spectral solvers for nonlocal equations on bounded domains. These spectral solvers exploit the fact that integration in the nonlocal formulation transforms into multiplication in Fourier space and that nonlocality is decoupled from the grid size, allowing fast and accurate solutions to the nonlocal problems. Our approach extends the spectral solvers developed in~\cite{alali2020fourier} for periodic domains by incorporating the two-dimensional Fourier continuation (2D-FC) algorithm introduced in~\cite{bruno2022two}. We evaluate the performance of the proposed methods on two-dimensional nonlocal Poisson and nonlocal diffusion equations defined on bounded domains. 
While the regularity of solutions to these equations in bounded settings remains an open  problem, we conduct numerical experiments to explore this issue, particularly focusing on studying discontinuities.
\end{abstract}


\section{Introduction and overview}
In this work, we present a spectral numerical method for solving two-dimensional nonlocal equations on bounded domains. In particular, we solve nonlocal Poisson and nonlocal diffusion equations within bounded domains. We consider nonlocal Laplace operators of the form
\begin{equation*}
	\Ldel u(x) = \cdel \int_{B_\delta(x)} \frac{u(y)-u(x)}{\|y-x\|^\beta}dy,~x\in\Rn,
\end{equation*}
where $\delta>0$ is called the horizon, $\beta$ satisfies $\beta<n+2$, and $\cdel$ is a scaling constant.

Nonlocal operators of this form have their roots in peridynamics \cite{bobaru2016handbook, silling2000reformulation, silling2007peridynamic, silling2017modeling} and they have been used in different applied settings, such as fractures  \cite{hu2012peridynamic, bobaru2015cracks, dias2017review, lipton2014dynamic}, image processing \cite{gilboa2009nonlocal}, tumor growth \cite{lejeune2017modeling}, material failure and damage \cite{isiet2021review, diehl2015simulation, hu2017peridynamics}, among others. There have been several mathematical and computational results involving nonlocal equations, see for example \cite{de2017finite, foss2016differentiability, foss2018existence,hinds2012dirichlet,mengesha2014nonlocal,mengesha2016characterization,burkovska2024efficient,madenci2014coupling}. On the mathematical side, the work in \cite{alali2021fourier} introduces the Fourier multipliers for nonlocal Laplace operators, studies the asymptotic behavior of these multipliers, and then applies the asymptotic analysis in the periodic setting to prove regularity results for the nonlocal Poisson equation. The Fourier multipliers $\mdel$ of the nonlocal operator $\Ldel$ are defined through Fourier transform by
\begin{equation}
	\widehat{\Ldel u} = \mdel \hat{u}.\label{multipliers_property}
\end{equation}

These multipliers are given by the integral representation~\eqref{multiplier_integral_rep}, and they are used to obtain the eigenvalues of the nonlocal operator $\Ldel$ in the case of periodic domains by evaluating them on a lattice of points (see Section 4.1 of~\cite{alali2021fourier}).
\begin{equation}
  \mdel(\nu)=\cdel\int_{B_\delta(0)}\frac{\cos(\nu\cdot z)-1}{\|z\|^\beta}dz. \label{multiplier_integral_rep}
\end{equation}
As pointed out in \cite{du2016asymptotically, du2017fast}, it is challenging to accurately and efficiently compute the eigenvalues of nonlocal operators. The work in \cite{alali2020fourier} simplifies the computation of these multipliers by representing the $n$-dimensional integral in \eqref{multiplier_integral_rep} by a hypergeometric function given by
\begin{equation}
    \mdel(\nu) = -\|\nu\|^2 ~_2F_3\left(1,\frac{n+2-\beta}{2};2,\frac{n+2}{2},\frac{n+4-\beta}{2};-\frac{1}{4}\|\nu\|^2\delta^2\right),
    \label{eq:multiplier-general}
\end{equation}
for any spacial dimension $n$ and any $\beta<n+2$. This reduces the computation of the multipliers in $\Rn$ to computing the 1D smooth function in \eqref{eq:multiplier-general}.

The work in \cite{mustapha2023regularity} and \cite{dang2024regularity} apply these Fourier multipliers to study the regularity of solutions for nonlocal diffusion equations and nonlocal wave equations over the space of periodic distributions, respectively. As emphasized in~\cite{mustapha2023regularity}, nonlocal diffusion equation satisfies an  instantaneous smoothing effect when the integral kernel is singular with  $\beta>n$ and a gradual (over time) smoothing effect when $\beta=n$. However, for  integrable kernels ($\beta<n$), the nonlocal diffusion equation is non-smoothing. According to the work in~\cite{alali2021fourier}, the regularity of the solution $u$ for the nonlocal Poisson equation $\Ldel u = f$ is the same as the regularity of $f$ when $\beta\le n$, while when $\beta>n$, the solution is $\beta-n$ more regular than $f$. Specifically, if $f\in H^s(T^n)$ for $s\ge 0$, then the solution $u\in H^s(T^n)$ when $\beta\le n$ and $u\in H^{s+\beta-n}(T^n)$ when $\beta>n$.

For computational results, on the other hand, many computational techniques and numerical analysis have been developed for solving nonlocal equations, see for example ~\cite{burkovska2024efficient, de2017finite, du2017fast, albin2011spectral, burkovska2020affine, du2016asymptotically, alali2020fourier, slevinsky2018spectral}. Spectral methods are accurate and efficient for solving nonlocal equations because they transform the integral in the nonlocal formulation to multiplication in Fourier space. For example, the application of the operator $\Ldel$ on a function $u$ results in multiplication in Fourier space, as shown in~\eqref{multipliers_property}, and this runs at a computational cost of $\mathcal{O}(N\log N)$ for an $N$-point discretization grid.
For comparison of the efficiency and accuracy of spectral methods versus standard numerical methods such as finite differences, see \cite{alali2020fourier}.

Some spectral methods for nonlocal equations have been developed in \cite{du2016asymptotically, alali2020fourier, jafarzadeh2020efficient,du2017fast,slevinsky2018spectral}. The work in \cite{jafarzadeh2020efficient} introduces a boundary-adapted spectral method for peridynamic diffusion problems with arbitrary boundary conditions. A spectral method for nonlocal equations on periodic domains has also been developed in \cite{alali2020fourier} with high accuracy. One difficulty that arises in adapting the spectral method developed in \cite{alali2020fourier} to solve nonlocal equations on bounded domains is that the solutions are in general nonperiodic, therefore exhibiting the Gibbs phenomenon as the Fourier coefficients of the solution decay slowly. One way around this is to use Fourier continuation. The idea behind Fourier continuation is to reduce the impact of slowly decaying Fourier coefficients of a smooth but nonperiodic function on a given domain by constructing a smooth periodic extension of the function over a larger domain.

Our proposed solvers build on the Fourier spectral method developed in \cite{alali2020fourier} for periodic domains and incorporate the two-dimensional Fourier continuation (2D-FC) algorithm from \cite{bruno2022two} to solve two-dimensional nonlocal equations on bounded domains. We also use this approach to perform numerical experiments aimed at exploring the regularity of solutions for nonlocal equations, with particular focus on  discontinuities. For integrable kernels, numerical simulations indicate that discontinuities may develop in the solution. In contrast, for singular kernels with $\beta>n+1=3$, the solution appear to remain continuous in the interior of the domain, and no discontinuities are observed.

The paper is organized as follows. In Section~\ref{2D-FC}, we  describe the two-dimensional Fourier continuation (2D-FC) algorithm developed in~\cite{bruno2022two}, along with minor modifications. Section~\ref{applications} presents our spectral solvers for two-dimensional nonlocal Poisson and nonlocal diffusion equations on bounded domains. In Section~\ref{numerical_results}, we discuss the numerical results and explore the behavior of solutions in various regimes. In Section~\ref{discussion}, we provide a summary of our approach and outline directions for future research.

\section{Two dimensional Fourier continuation}\label{2D-FC}
In this section, we provide a description of the 2D-FC algorithm presented in \cite{bruno2022two}. The description includes a modification of some of the steps in the algorithm to improve efficiency and accuracy. In particular, we modify the steps in computing the continuation values on the Cartesian grids in the extended region. For consistency, we follow a notation similar to that of~\cite{bruno2022two}.

\subsection{Overview of 1D blending-to-zero algorithm} \label{blend-to-zero}
\begin{figure}[b!]
    \centering
    \begin{subfigure}[h]{0.92\textwidth}
        \centering
        \includegraphics[width=\textwidth]{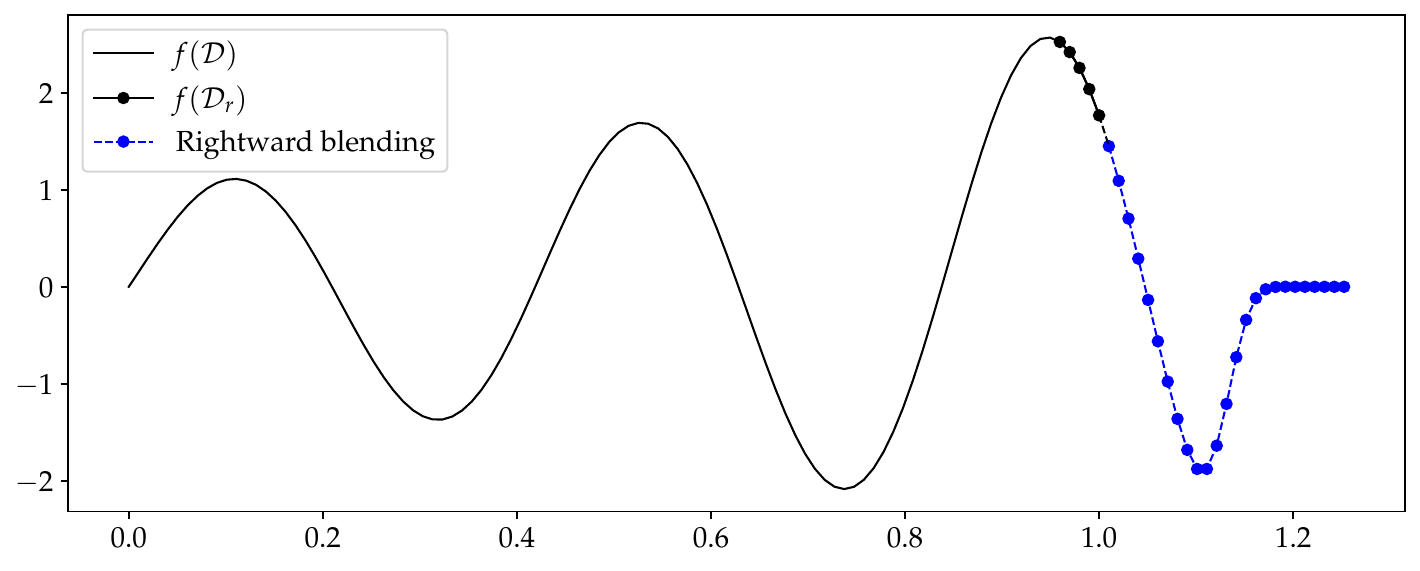}
    \end{subfigure}
    \caption{A demonstration of the 1D blending to zero on a smooth but nonperiodic function $f(x)=\sin(15x)\exp(x)$ defined on [0,1]. The black solid points denote the $d$ near-boundary function values used to produce the extension of $f$ denoted by the blue solid points.}
    \label{blendtozero}
\end{figure}

In the description given in \cite{albin2014discrete}, the blending-to-zero algorithm starts by smoothly extending the vector $f=(f_1,f_2,\dots,f_N)\in\bbR^N$ of values of a complex-valued function $f$ over the set of uniformly spaced nodes  $\mathcal{D} = \{x_1,x_2,\dots,x_N\}$ where for a given step size $h>0$, $x_j = x_1+(j-1)h.$
Without loss of generality, we assume that $x_1=0$ and $x_N=1$. The extension relies on some number $d$ of near-boundary function values, $\{f_1,f_2,\dots,f_d\}$ near $x=0$ and $\{f_{N-d+1},f_{N-d+2},\dots,f_N\}$ near $x=1$, and it produces a sequence $\{\tilde{f}_j\}_{j=-\infty}^\infty$ with the property that 
\[\tilde{f}_j = f_j \text{ for } j=1,2,\dots,N.\] 
The extended vector $\tilde{f}$ is then blended to zero on either end by multiplying with a suitable windowing function. A complete description of this algorithm is provided in \cite{albin2014discrete}. In Figure~\ref{blendtozero}, we demonstrate the 1D blending to zero of the nonperiodic function $f:[0,1]\to\bbR$ given by $f(x)=\sin(15x)e^x$.

\subsection{2D-FC algorithm}
\begin{figure}[t!]
    \centering
    \begin{subfigure}[h]{0.4\textwidth}
        \centering
        \includegraphics[width=\textwidth]{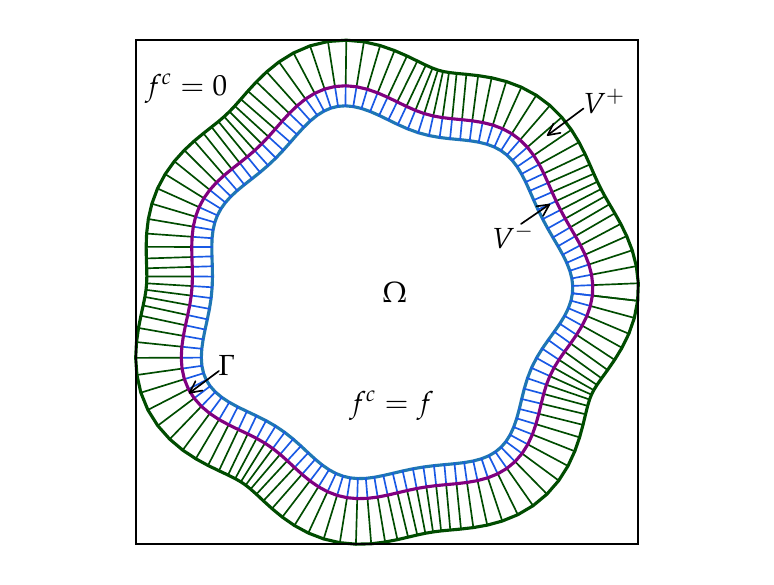}
    \end{subfigure}
    \caption{Geometric constructions of the 2D-FC procedure. The interior strip $V^-$ and the exterior strip $V^+$ are defined in subsection~\ref{curvilinear_grid}.}
    \label{2D_algo}
\end{figure}

Let $f:\overline{\Omega}\to \mathbb{C}$ be given, where the domain $\Omega\subset \bbR^2$ has a smooth boundary $\Gamma=\overline{\Omega}\setminus \Omega$. We assume that values of $f$ are known at the uniform Cartesian grid points within $\overline{\Omega}$ and at the grid points on the boundary. The algorithm first produces 1D blending to zero values for $f$ along directions normal to the boundary $\Gamma$, to produce continuation values that vanish on a certain 2D tangential-normal curvilinear grid around $\Gamma$ as depicted in Figure~\ref{2D_algo}. These continuation values are then interpolated onto Cartesian grid points inside the region $V^+$ to produce a 2D blending to zero continuation of $f$. Once all the interpolated values have been obtained at all Cartesian grid points around $\Gamma$, the 2D-FC function is then obtained via a 2D Fast Fourier Transform (FFT) given by
\begin{eqnarray}
	f^c(x,y) = \sum_{k = -\frac{N_x}{2}+1}^{N_x/2}\sum_{l = -\frac{N_y}{2}+1}^{N_y/2} \hat{f}^c_{kl}\exp\left( \frac{2\pi i kx}{L_x}\right)\exp\left( \frac{2\pi i ly}{L_y}\right), \label{fourier_series}
\end{eqnarray}
where $L_x$ and $L_y$ are the periods in $x$ and $y$ directions respectively. The steps of this algorithm are described in the following subsections.

\subsubsection{Curvilinear grid along a 2D tangential-normal}\label{curvilinear_grid}
The curvilinear grids are introduced within the interior and exterior strips $V^-$ and $V^+$ respectively, as shown in Figure~\ref{2D_algo}. These strips are given by 
\begin{eqnarray*}
    V^- &=& \{q(\theta)-n(\theta)\gamma: 0\le \theta\le 2\pi \text{ and } 0\le \gamma \le (d-1)k_1\} \\
    V^+ &=& \{q(\theta)+n(\theta)\gamma: 0\le \theta\le 2\pi \text{ and } 0\le \gamma \le Ck_1\},
\end{eqnarray*}
where $q(\theta) = (x(\theta), y(\theta))$ is the parametrization of $\Gamma$, $n(\theta) = (n_x(\theta),n_y(\theta))$ is the unit normal to the boundary at $q(\theta)$, and $d,C,k_1$ represent the number of matching points near the boundary, the number of continuation points, and the step size along the normals respectively. Consider the following uniform discretization of the interval $[0,2\pi]$
\begin{eqnarray}
    \Theta_B = \{\theta_p=pk_2: 0\le p < B, ~k_2 = 2\pi/B\}. \label{bdry_discretization}
\end{eqnarray}
Then for each $\theta_p\in\Theta_B$, the curvilinear grid along the corresponding normal within $V^-$ and $V^+$ are respectively given as
\begin{eqnarray*}
    V^-_{\theta_p}&=&\{r_{p,q}=q(\theta_p)+n(\theta_p)(q-d+1)k_1, ~0\le q\le d-1\}\\
    V^+_{\theta_p}&=&\{s_{p,q} = q(\theta_p)+n(\theta_p)qk_1/n_r, ~0\le q\le C_r\},
\end{eqnarray*}
where $C_r = Cn_r$ for some integer refinement factor $n_r$ (Remark~\ref{refine_factor}). Let $\mathcal{R} = [a_0,a_1]\times [b_0,b_1]$ be the smallest closed rectangle which contains $\Omega\cup V^+$ as shown in Figure~\ref{2D_algo}, and define the following Cartesian grid $\mathcal{H}$ on $\mathcal{R}$ with step size $h$, where the continuation function values will be computed.
\begin{eqnarray*}
    \mathcal{H}=\{z_{ij}=(x_i,y_j): x_i=a_0+ih; ~y_j = b_0+jh, 0\le i<N_x \text{ and }0\le j<N_y\}.
\end{eqnarray*}
\begin{remark}\label{refine_factor}
	The refinement factor $n_r$ helps to prevent accuracy loss by making the normal-direction grids at the exterior finer than the grids inherent by the blending-to-zero process.
\end{remark}

\subsubsection{Computing continuation values on $V^+_{\theta_p}$}\label{FC_values}
\begin{figure}[h!]
    \centering
    \begin{subfigure}[h]{0.3\textwidth}
        \centering
        \includegraphics[width=\textwidth]{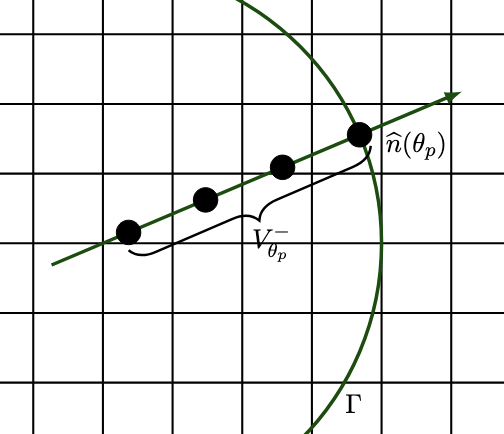}
    \end{subfigure}
    \hspace{20pt}
    \begin{subfigure}[h]{0.3\textwidth}
        \centering
        \includegraphics[width=\textwidth]{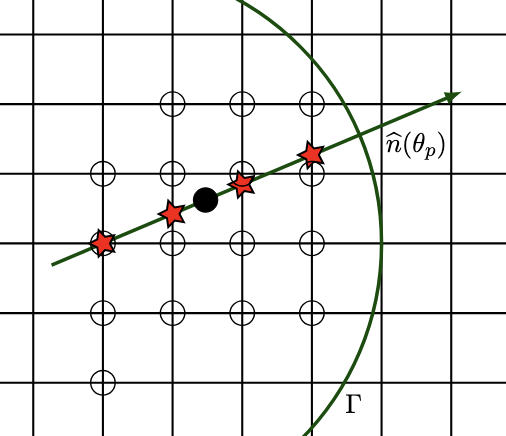}
    \end{subfigure}
    \caption{Interpolation procedure for evaluating values of $f$ on the set $V_{\theta_p}^-$ for the case $|n_x(\theta_p)|\ge|n_y(\theta_p)|$. The black solid circles indicate the points that define the set $V^-_{\theta_p}$. The known values of $f$ at the $M$ open circles are used to interpolate the function values at the red stars (intersection points of the normal and the vertical grid lines). The values of $f$ at the red points are then used to obtain the values of $f$ at the black points.}
    \label{continuation_along_normals}
\end{figure}

For a given $\theta_p\in \Theta_B$, a continuation of $f$ on $V^+_{\theta_p}$ is obtained via the blending-to-zero algorithm described in Section~\ref{blend-to-zero}. We first obtain the approximations of $f$ on the set $V^-_{\theta_p}$, which are indicated by black solid circles in Figure~\ref{continuation_along_normals}, and then use these approximations to obtain a continuation of $f$ via the 1D blending-to-zero procedure along the corresponding normal direction. The values of $f$ on $V^-_{\theta_p}$ are obtained by a two-step polynomial interpolation scheme, using polynomials of degree $M-1$. This is briefly described as follows:

Consider the case where $|n_x(\theta_p)|\ge |n_y(\theta_p)|$. We first estimate the values of $f$ at each red point, which is the intersection of the normal line and the vertical Cartesian grid lines as shown in Figure~\ref{continuation_along_normals}, by interpolating vertically the function values at the selected $M$ Cartesian points, shown as open circles, onto the red points. Then we interpolate the values of $f$ at the $M$ red points onto the black solid circles, using a polynomial of degree $M-1$. 

The case where $|n_x(\theta_p)| < |n_y(\theta_p)|$ is handled in a similar way by replacing the interpolation along vertical Cartesian grid lines with interpolation along horizontal Cartesian grid lines to obtain the function values at the red points, and then interpolating the function values at the red points onto the black solid points. 

\subsubsection{Computing continuation values on the Cartesian grid}  \label{FC_on_cartpts}
Once all the continuation values have been obtained on $V_{\Theta_B}^+ = \{V_{\theta_p}^+: \theta_p\in \Theta_B\}$,
we then obtain the continuation of $f$ on the Cartesian grid points within $V^+$ using a two-step 1D polynomial interpolation scheme. The 2D-FC algorithm is then completed by evaluating corresponding FCs in \eqref{fourier_series} by means of 2D FFT. The description of the two-step 1D polynomial interpolation scheme is presented below.

For an arbitrary Cartesian point $Q\in \mathcal{H}\cap V^+$, we first obtain the nearest boundary parameter value, $\theta_p$, around $Q$ using a proximity map. This proximity map can be obtained by associating with each curvilinear discretization point $s_{p,q}$ the nearest Cartesian point using \textit{round} operator (which rounds $x$ and $y$ coordinates of $Q$ to the nearest integer), resulting in a set $P_0\subseteq (\mathcal{H}\cap V^+)\times V_{\Theta_B}^+$ of pairs of points. The set $P_0$ is then modified by removing multiple associations for a given Cartesian point, and, by adding a next nearest curvilinear neighbor to the Cartesian points that are unassociated. The algorithm then uses the $M$-nearest boundary parameter values
\begin{equation*}
	S_{\theta_p} = \{\theta_{p-K_l}, \theta_{p-K_l+1},\dots,\theta_p,\theta_{p+1},\dots,\theta_{p+K_r}\}, 
\end{equation*}
around $\theta_p$ (where $K_l+K_r+1 = M$) to obtain the continuation of $f$ at $Q$. Parameter values $\theta_k$ with negative values of $k$ are interpreted by periodicity: $\theta_k = \theta_{B+k}$.

\begin{figure}[b!]
    \centering
    \begin{subfigure}[h]{0.37\textwidth}
        \centering
        \includegraphics[width=\textwidth]{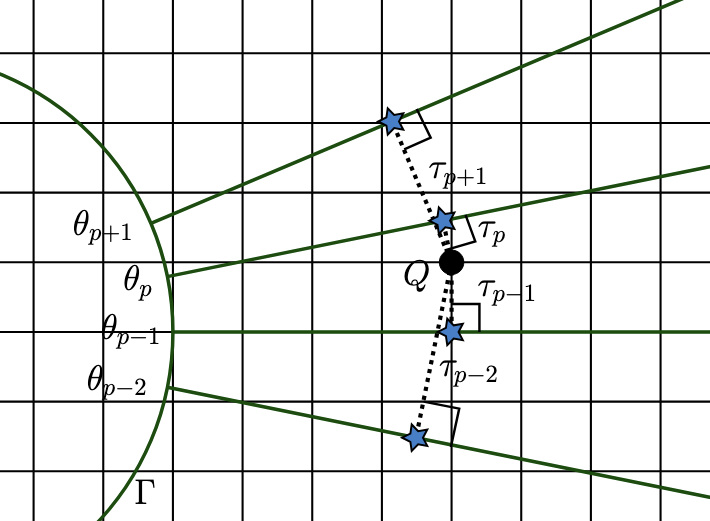}
        \caption{}
        \label{bluepoints1}
    \end{subfigure}
    \hspace{20pt}
    \begin{subfigure}[h]{0.37\textwidth}
        \centering
        \includegraphics[width=\textwidth]{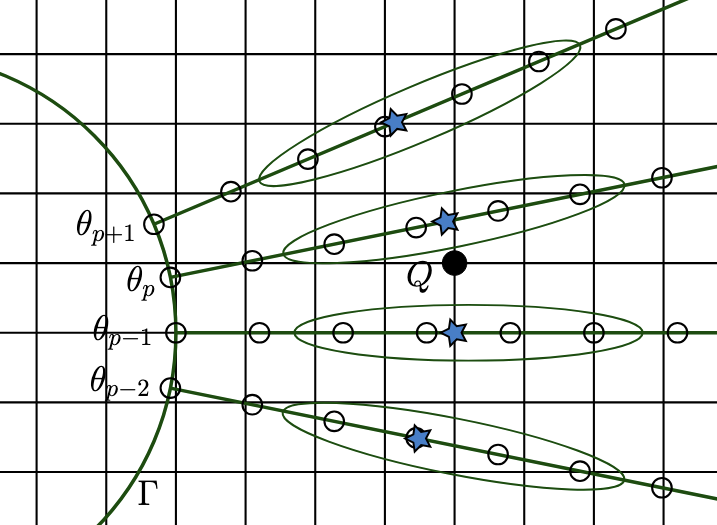}
        \caption{}
        \label{bluepoints2}
    \end{subfigure}
    \caption{Interpolation procedure for obtaining continuation values of $f$ on $\mathcal{H}\cap V^+$. In (a), the blue-star interpolation points are the perpendicular projections of the Cartesian point $Q$ onto the normal vectors, and the $\tau$'s are the signed perpendicular distances from $Q$ to each normal vector. In (b), the interpolation of continuation values from the curvilinear grids to a point $Q$ on the Cartesian~grids.}
    \label{bluepoints}
\end{figure}

For each parameter value $\theta_j, ~p-K_l\le j\le p+K_r$, we find the signed perpendicular distance from point $Q$ to the normal vector $n(\theta_j)$, denoted by $\tau_j$ in Figure~\ref{bluepoints1}. Let $b_j$ be the perpendicular projection of $Q$ onto the normal vector $n(\theta_j)$, denoted as the blue points in Figure~\ref{bluepoints}. Define
\begin{equation*}
	\mathcal{T}(\tau_j) = f^c(b_j),~p-K_l\le j\le p+K_r,
\end{equation*}
where $f^c(b_j)$ is the continuation of $f$ at $b_j$ obtained via a 1D polynomial interpolation on $M$  blending to zero points along $n(\theta_j)$ (see Figure~\ref{bluepoints2}). Clearly, $\mathcal{T}(0) = f^c(Q)$. It then follows that a 1D polynomial interpolation procedure on the function $\mathcal{T}(\tau)$, using polynomial of degree $M-1$, can be used to obtain the desired continuation of $f$ at Cartesian point $Q$.  

Finally, by applying a 2D FFT to the continuation function values computed above, we obtain the desired Fourier series expression in \eqref{fourier_series} for the continuation function.

\begin{remark}\label{higher_normals}
	We observe that when the value of $M$ becomes too large, the accuracy decreases. This loss of accuracy arises because using a larger number of normals in the interpolation  procedure for estimating the point $Q$ (see Figure~\ref{bluepoints}) causes the corresponding  blue points to move closer to the boundary of $\Omega$. Consequently, the estimation of $Q$ is affected, since the values of the blue points near the boundary are not sufficiently blended to zero. To avoid this, we select a smaller number of normals (corresponding to a lower $M$) so that the contributing normals remain close to the point~$Q$.
\end{remark}

\subsubsection{Parameter description and selection}\label{parameter_selction}
Here, we briefly describe the parameters used in the 2D-FC algorithm. Additionally, we select some specific values for these parameters to get the best result.
\begin{itemize}[noitemsep]
	\item $d$: number of near-boundary points used to blend $f$ to zero along each normal. 
	\item $C$: number of continuation points in the 1D blending to zero process. 
	\item $k_1$: normal step size. 
	\item $k_2$: boundary step size. 
	\item $h$: step size used to construct Cartesian grid $\mathcal{H}$. 
	\item $n_r$: refinement factor used to obtain the curvilinear grids along each normal vector (Remark~\ref{refine_factor})
	\item $N=N_x\times N_y$: number of points in the uniform Cartesian grid 
	\item $B$: number of boundary discretization points. 
	\item $M-1$: the degree of the interpolating polynomial. 
\end{itemize}

For a given step size $h$ used to construct the Cartesian grid $\mathcal{H}$, we choose normal step size $k_1$ and boundary step size $k_2$ to coincide with $h$. The parameters $C=25$ and $n_r=6$ were used in the 1D blending to zero procedure described in Subsection~\ref{blend-to-zero}. Finally, we set the interpolating polynomial degree $(M-1)=d$ for all matching point numbers $d$ considered except the interpolation-degree experiments presented in Example~\ref{performance}.

In the following example, we demonstrate the performance of the 2D-FC method using three interpolation degrees: $M-1=d,$ $d+1$, and $d+2$ for a given matching point number $d$.

\begin{ex}\label{performance}
	We consider the FC approximation of a function $f:\Omega\to \bbR$ given by
	\begin{equation}
		f(x,y)=-(x^8+y^8)\sin(8\pi x)\sin(8\pi y)\label{f_performance}
	\end{equation}
	on the unit disk $\Omega=\{(x,y):x^2+y^2\le 1\}$ with three different values of $d$: $d=4$, $5$, and $8$; and three different values of $M$: $M=d+1$, $d+2$, and $d+3$.
\end{ex}

\begin{figure}[t!]
    \centering
    \begin{subfigure}[h]{0.325\textwidth}
        \centering
        \includegraphics[width=\textwidth]{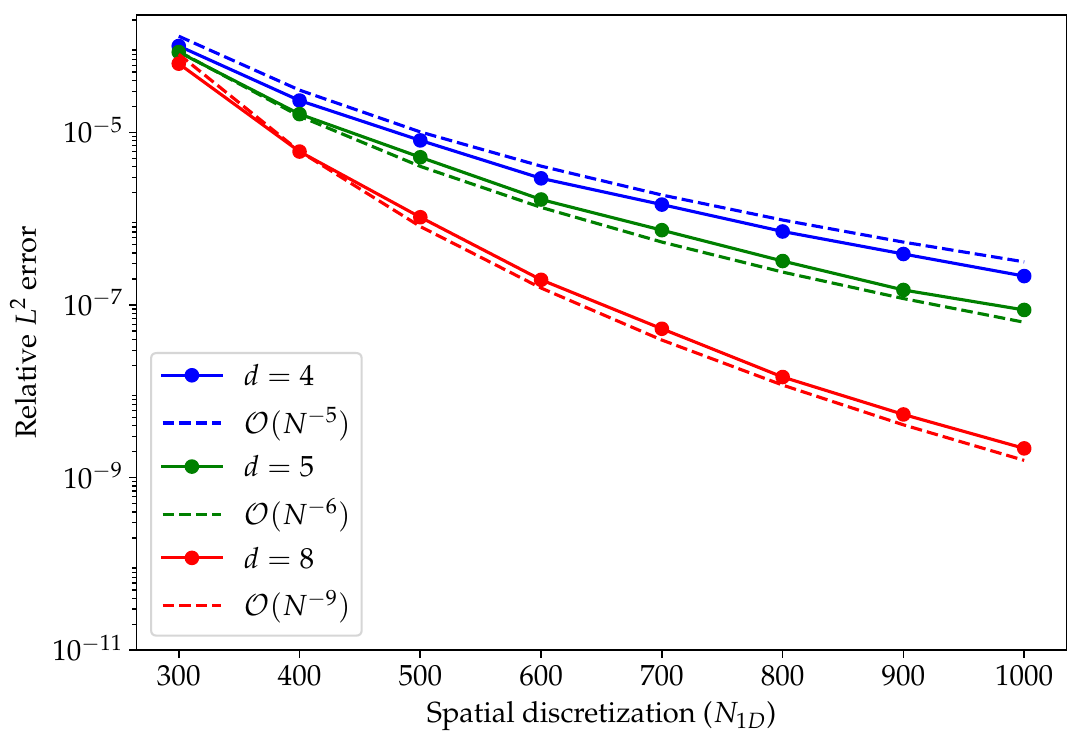}
        \caption{$M=d+1$}
        \label{dplus1}
    \end{subfigure}
    \begin{subfigure}[h]{0.325\textwidth}
        \centering
        \includegraphics[width=\textwidth]{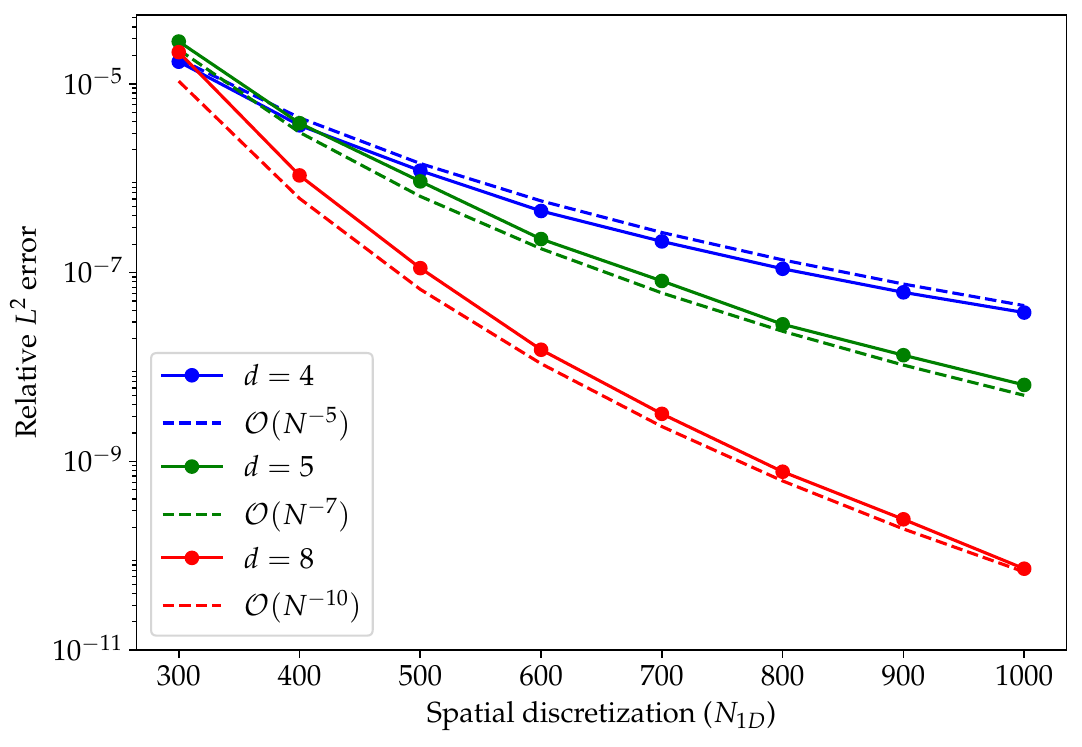}
        \caption{$M=d+2$}
        \label{dplus2}
    \end{subfigure}
    \begin{subfigure}[h]{0.325\textwidth}
        \centering
        \includegraphics[width=\textwidth]{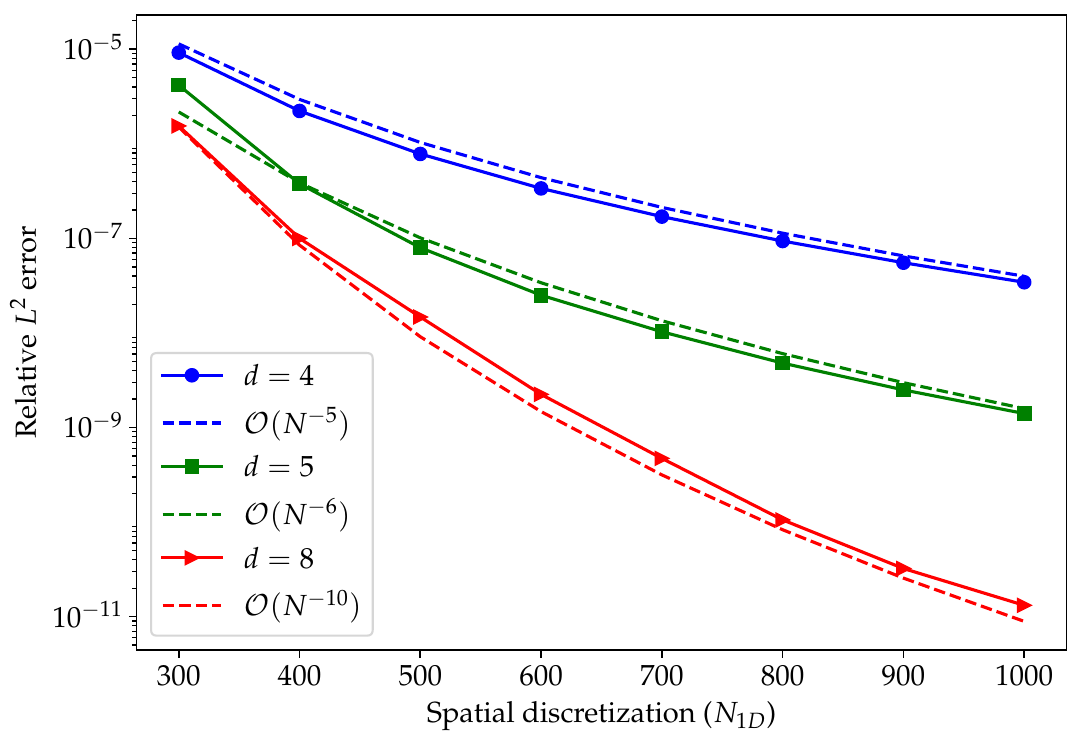}
        \caption{$M=d+3$}
        \label{dplus3}
    \end{subfigure}
    \caption{Convergence study for the 2D-FC approximation of the function $f$ in \eqref{f_performance} with three choices of matching point numbers $d$: $d=4$, $5$, and $8$; and three choices of $M$: $M=d+1$, $d+2$, and $d+3$. The integer $N_{1D}$ is the number of spatial grid points in one spatial direction within $\Omega$.}
    \label{performance_2DFC}
\end{figure}

Figure~\ref{performance_2DFC} displays the relative $L^2$ errors $\mathcal{E}_{L^2}^{\text{rel}}$, given in~\eqref{l2error},  as functions of the number of spatial discretization, $N_{1D}$ ($N_x=N_y\approx N_{1D}+2C$), in one direction within the domain $\Omega=\{(x,y):x^2+y^2\le 1\}$.
The errors are computed on a Cartesian grid of step size $h/2$ within $\Omega$. The approximations of $f$ on the finer grid are obtained by means of zero padding~\cite{press2007numerical}. Specifically, we compute $\mathcal{E}_{L^2}^{\text{rel}}$ using the following formula:
\begin{equation}
	\mathcal{E}_{L^2}^{\text{rel}} = \frac{\|f-f^{\mathrm{zp}}\|_{L^2(\Omega)}}{\|f\|_{L^2(\Omega)}},\label{l2error}
\end{equation}
where $f^{\mathrm{zp}}$ denotes the zero--padded Fourier approximation of $f$, obtained by extending 
$f$ through Fourier continuation, then applying zero padding in Fourier space, and restricting the result back to $\Omega$. 

For the case $M=d+1$ (Figure~\ref{dplus1}), we observe convergence of orders five, six, and nine for $d=4$, $d=5$, and $d=8$ respectively, corresponding to an overall convergence rate of $\mathcal{O}(h^{d+1})$. However, for the case $M=d+2$ (Figure~\ref{dplus2}), the observed rates are $\mathcal{O}(h^{d+1})$ for $d=4$ and $\mathcal{O}(h^{d+2})$ for $d=5$ and $d=8$. Similarly, for $M=d+3$ (Figure~\ref{dplus3}), we observe $\mathcal{O}(h^{d+1})$ convergence for $d=4$ and $d=5$ and $\mathcal{O}(h^{d+2})$ convergence for $d=8$. In particular, we observe a lower convergence rate for higher values of $M$ as explained in Remark~\ref{higher_normals}.
The choice $M=d+1$, which yields a consistent convergence order of $d+1$ for the three values of $d$ considered, will be employed in all subsequent numerical experiments.

In the next example, we demonstrate the 2D-FC expansion of a two-dimensional nonperiodic function defined on a two-dimensional domain contained within the curve given by
\begin{equation}
	x(\theta)=r(\theta)\cos(\theta);~y(\theta)=r(\theta)\sin(\theta)
	\label{ringdomain},
\end{equation}
where $r(\theta) = 5 +\cos(7\theta)/2 + \sin(4\theta)/3$ and $0\le \theta \le 2\pi.$

\begin{figure}[ht!]
    \centering
    \begin{subfigure}[h]{0.42\textwidth}
        \centering
        \includegraphics[width=\textwidth, height=4.2cm]{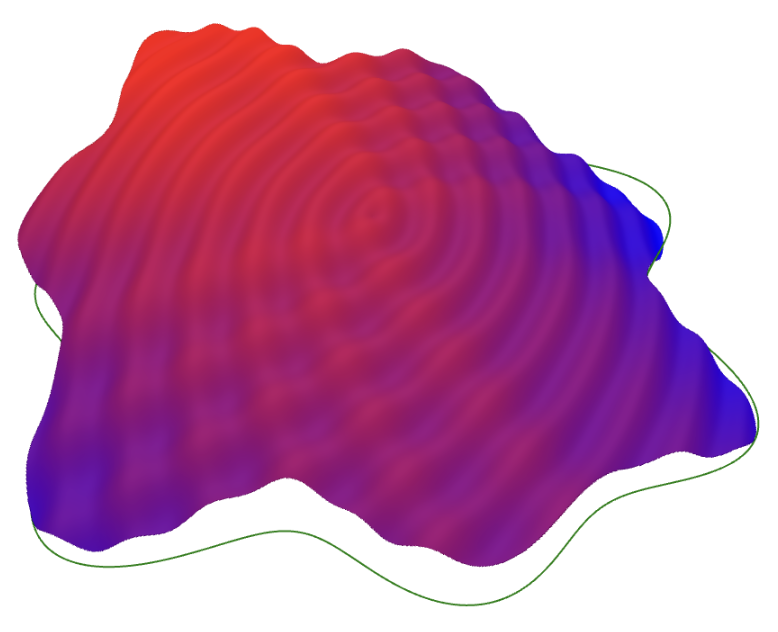}
        \caption{}
    \end{subfigure}
    \hfill
    \begin{subfigure}[h]{0.5\textwidth}
        \centering
        \includegraphics[width=\textwidth, height=4.2cm]{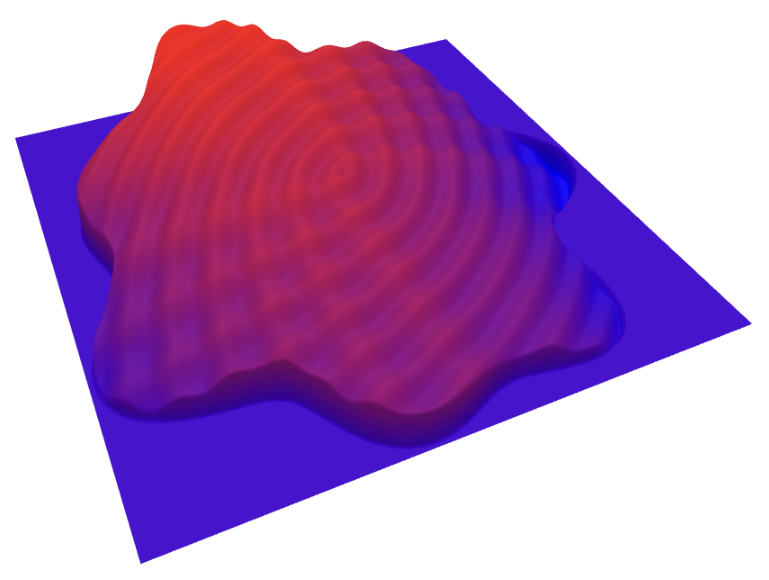}
        \caption{}
    \end{subfigure}
    \caption{Demonstration of the 2D-FC expansion of a nonperiodic function defined on a 2D domain in the settings of Example~\ref{ex:2D-FC-demo}. (a) the original nonperiodic function $g$. The green curve indicates the boundary of the domain. (b) the 2D-FC expansion of the function $g$.}
    \label{2D-FC-demo}
\end{figure}

\begin{ex}[Graphical illustration of 2D-FC method]\label{ex:2D-FC-demo}
	We choose $d=4$, $M=5$, $h=0.01$, and employ the 2D-FC method to extend the function
	\begin{equation*}
		g(x,y) = 6 + 0.5y - 0.2(x+1)^2 + 0.6\sin(4.1\sqrt{x^2+y^2})\cos(3.8(x-y)) ,
	\end{equation*}
	over the domain defined in \eqref{ringdomain}. 
\end{ex}
Figure~\ref{2D-FC-demo} shows the original function together with its periodic extension obtained using the 2D-FC method. The corresponding relative $L^2$ error $\mathcal{E}_{L^2}^{\text{rel}}$, computed using the formula in~\eqref{l2error}, is $1.1674\cdot 10^{-06}$.

\section{Applications of 2D-FC methods}\label{applications}
\begin{figure}[ht!]
	\centering
	\begin{subfigure}[h]{0.33\textwidth}
		\centering
		\includegraphics[width=\textwidth]{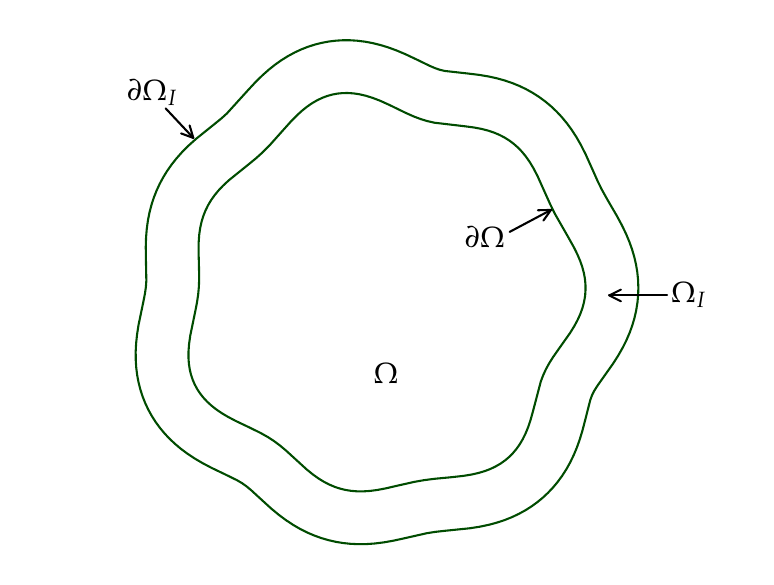}
	\end{subfigure}
	\caption{The collar region}
	\label{collar_region}
\end{figure}

In this section, we present a 2D-FC based method for solving two-dimensional nonlocal Poisson and nonlocal diffusion equations within a bounded domain  $\Omega$ and with Dirichlet conditions on the collar region. The collar region is the extension of $\Omega$ by $\delta$ as shown in Figure~\ref{collar_region} and is defined by
\begin{equation*}
	\collar = \{{\bf{y}}\in\bbR^2\setminus\Omega:\|{\bf{x}}-{\bf{y}}\|\le \delta, \text{ for some }{\bf{x}}\in\Omega\}.
\end{equation*}
A description of our proposed method is given in the following sections.
\subsection{Nonlocal Poisson equation}\label{sec:nonlocal_poisson}
Consider the following 2D nonlocal Poisson equation with Dirichlet conditions on the collar region
\begin{eqnarray}
	\begin{cases}
		\Ldel u(x,y) = f(x,y),& (x,y)\in \Omega \\
		u(x,y) = b(x,y), & (x,y)\in \collar
	\end{cases}\label{nonlocal_poisson}.
\end{eqnarray}
Our approach for solving \eqref{nonlocal_poisson} builds on the spectral method for nonlocal Poisson equation on periodic domains by constructing a periodic extension of the problem using the Fourier continuation algorithm. We then solve the extended problem using the iterative linear algebra solver GMRES~\cite{saad2003iterative} to obtain the solution within a specified tolerance. This process is detailed below.

We first transform \eqref{nonlocal_poisson} into the following linear system
\begin{eqnarray}
	Au = g, \quad A = \begin{pmatrix}L \\ S\end{pmatrix}, \quad g = \begin{pmatrix}f \\ b\end{pmatrix},\label{associated_system_poisson}
\end{eqnarray}
where $L$ and $S$ are respectively the associated matrices of the nonlocal operator $\Ldel$ and the boundary operator $\mathcal{S}$ given by the Dirichlet conditions on the collar region. To apply GMRES to this system, we do not need to compute $A$ and store it in memory. Instead, it is sufficient to express the forward operator $v \mapsto Av$ as a function. The action of the operator $\mathcal{S}$ is a restriction of the values of $v$ to the collar region. We derive the forward operator function using the properties of $\Ldel$ as follows:
\begin{equation}
	u\to u_{\text{ext}}\to \widehat{u}_\text{ext}\to \widehat{Lu_\text{ext}} = \mdel \widehat{u}_{\text{ext}}\to Lu_\text{ext}\to Lu,\label{eq:FC-Poisson-steps}
\end{equation}
where $u_{\text{ext}}$ is the extension of $u$ obtained by the 2D-FC method, $\widehat{u}_\text{ext}$ is the Fourier transform of $u_\text{ext}$ obtained by means of the fast Fourier transform (FFT), $Lu_\text{ext}$ is the inverse Fourier transform of $\widehat{Lu_\text{ext}}$, and $Lu$ is the restriction of $Lu_\text{ext}$ to the domain $\Omega$. The extension step of~\eqref{eq:FC-Poisson-steps} has a computational complexity of $\mathcal{O}(N^{\frac{1}{2}})$; each interpolation and each blending-to-zero has constant complexity. Applying the operator in Fourier space is a $\mathcal{O}(N)$ operation, and the final restriction operator is $\mathcal{O}(1)$. Thus, the total computational complexity of~\eqref{eq:FC-Poisson-steps} is dominated by the $\mathcal{O}(N\log N)$ FFT. We obtain the solution of \eqref{nonlocal_poisson} by applying GMRES to the system~\eqref{associated_system_poisson} within a specified tolerance.

\subsection{Nonlocal diffusion equation}\label{sec:nonlocal_diffusion}
Consider the following 2D nonlocal diffusion equation with Dirichlet conditions on the collar region and with initial condition
\begin{eqnarray}
	\begin{cases}
		u_t(x,y,t) = \Ldel u(x,y,t) + s(x,y,t),& (x,y,t)\in \Omega\times (0,T), \\
		u(x,y,0) = u_0(x,y),&(x,y)\in \Omega,\\
		u(x,y,t) = b(x,y), & (x,y)\in \collar.
	\end{cases}\label{nonlocal_diffusion_eqn}
\end{eqnarray}
Similar to the nonlocal Poisson solver, we first construct the periodic extension of \eqref{nonlocal_diffusion_eqn} using the 2D-FC algorithm. We then solve the extended problem using the fourth-order Adams-Bashforth method. The process is detailed below.

First, consider the semi-discrete form of \eqref{nonlocal_diffusion_eqn}, where we replace $u$ by a time dependent vector $\Vec{u}$ as follows
\begin{eqnarray}
	\vec{u}'= \Ldel \vec{u}  + \vec{s} := f(\vec{u}).
	\label{diffusion-rhs}
\end{eqnarray}
For a given final time $T>0$ and a time step size $\tau$, we first obtain the vectors $\vec{u}_1$, $\vec{u}_2$, and $\vec{u}_3$, using the fourth-order Runge-Kutta method. Then, for $k=4,5,\dots, K$, where $K=T/\tau$, we obtain $\Vec{u}_k$ by using the fourth-order Adams-Bashforth method given by
\begin{equation*}
	\vec{u}_{k+4} = \vec{u}_{k+3} + \tau \left(\frac{55}{24}f(\vec{u}_{k+3})-\frac{59}{24}f(\vec{u}_{k+2})+\frac{37}{24}f(\vec{u}_{k+1})-\frac{9}{24}f(\vec{u}_{k})\right).
\end{equation*}
To obtain $f(\vec{u}_k)$, $k=0,1,\dots,K$, we apply the 2D-FC method to extend $\vec{u}$ in \eqref{diffusion-rhs}, as follows
\begin{equation*}
	\vec{u}_k \to (\vec{u}_k)_{\text{ext}}\to \widehat{(\vec{u}_k)}_{\text{ext}}\to \widehat{(\Ldel \vec{u}_k)}_{\text{ext}} = \mdel \widehat{(\vec{u}_k)}_{\text{ext}} \to \Ldel \vec{u}_k + \vec{s}_k = f(\vec{u}_k),
\end{equation*}
where $(\vec{u}_k)_{\text{ext}}$ is the extension of $\vec{u}_k$ obtained by the 2D-FC method, $\widehat{(\vec{u}_k)}_{\text{ext}}$ is the Fourier transform of $(\vec{u}_k)_{\text{ext}}$ obtained by means of FFT, $(\Ldel \vec{u}_k)_{\text{ext}}$ is the inverse Fourier transform of $\widehat{(\Ldel \vec{u}_k)}_{\text{ext}}$, and $\Ldel\vec{u}_k$ is the restriction of $(\Ldel \vec{u}_k)_{\text{ext}}$ to the domain $\Omega$.

\section{Numerical results}\label{numerical_results}
This section demonstrates the accuracy and performance of our proposed method for solving nonlocal equations on bounded domains. In particular, we present the convergence study for the spectral solvers described in subsections~\ref{sec:nonlocal_poisson}~and~\ref{sec:nonlocal_diffusion}  for the 2D nonlocal Poisson and nonlocal diffusion equations, respectively. Additionally, we study the discontinuities in the solutions of these equations. The 2D-FC parameters are as chosen in subsection~\ref{parameter_selction}.

\subsection{Accuracy and convergence test for nonlocal Poisson equation}
Here, we construct a manufactured solution to study the accuracy and convergence of the nonlocal Poisson equation~\eqref{nonlocal_poisson}. To this end, let 
\begin{equation*}
	u(x,y) = \sin(2\pi r_1x)\sin(2\pi r_2y)
\end{equation*}
be the solution of \eqref{nonlocal_poisson}, where $r_1= 10.6418$ and $r_2=12.6418$. Clearly, $u(x,y)$ is a highly oscillatory nonperiodic function. To obtain the right hand side function, $f(x,y)$, we apply the operator $\Ldel$ on $u$, using the property of $\Ldel$ given in \eqref{multipliers_property},
with fixed $\delta=0.4$ and three different values of $\beta:$ $1.2$, $2.0$, and $2.5$. These values of $\beta$ represents $\beta<n$, $\beta=n$, and $\beta>n$. A direct computation of the multipliers $\mdel$ is provided in Table~\ref{manufactured_sol_multipliers}.

\begin{table}[ht!]
	\centering
	\begin{tabularx}{\textwidth}{>{\centering\arraybackslash}c|>
			{\centering\arraybackslash}X|>
			{\centering\arraybackslash}X|>
			{\centering\arraybackslash}X}
		\hhline{|====|}
		&&&\\[-2ex]
		\hspace{7pt}$\delta$\hspace{7pt} & $\overset{\beta = 1.2}{\mdel}$ & $\overset{\beta=2.0}{\mdel}$ & $\overset{\beta=2.5}{\mdel}$ \\[0.5ex]
		\hline
		&&&\\[-2ex]
		$0.4$& $-82.87098585883194$ & $-180.5053934013443$ & $-387.0397711705603$\\
		[0.5ex]
		\hhline{|====|}
	\end{tabularx}
	\caption{The multipliers of the nonlocal operator $\Ldel$ applied on the function $u=\sin(2\pi r_1 x)\sin(2\pi r_2 y)$, where $r_1= 10.6418$ and $r_2 = 12.6418$, with a fixed $\delta=0.4$ and three different values of $\beta$.}
	\label{manufactured_sol_multipliers}
\end{table}
We compute these multipliers using the implementation of the $_2F_3$ hypergeometric function provided by mpmath \cite{johansson2013mpmath}, as detailed in \cite{alali2021fourier}. We then use these manufactured solutions to compute the errors as shown in the following example.

\begin{ex}[Convergence test for the nonlocal Poisson equation]\label{convergence_test_poisson}
	We employ our proposed spectral method to solve the nonlocal Poisson equation $\eqref{nonlocal_poisson}$ over the kite-shaped domain contained within the curve given by
	\begin{eqnarray}
		x(\theta) = \cos(\theta) + 0.35\cos(\theta)-0.35;~y(\theta) = 0.7\sin(\theta)~\text{ for } 0\le \theta\le 2\pi \label{kite_domain}
	\end{eqnarray} 
	and study the error convergence using step sizes $h=0.02/2^k,~k=0,1,2,3,4$. We fix $\delta=0.4$, and choose three different values of $\beta:$ $1.2$, $2.0$, $2.5$. We set the boundary function $b(x,y) = \sin(2\pi r_1x)\sin(2\pi r_2y)$ and the right-hand side function
	\begin{equation*}
		f(x,y) = \Ldel b(x,y) = \mdel b(x,y),
	\end{equation*}
	where $r_1=10.6418$, $r_2=12.6418$, and the multipliers $\mdel$ are as given in Table~\ref{manufactured_sol_multipliers}.
\end{ex}

A convergence study for this test case is shown in Table~\ref{conv_study_nonlocal_Poisson} for $d=4$ and $d=5$, with fixed $M=d+1$. We use the manufactured solution described above to compute the $L^2$ errors reported. We observe sixth and seventh orders of convergence for $d=4$ and $d=5$, respectively, showing convergence of order $d+2$ for the three cases of $\beta$ considered. 
Moreover, we observe no significant difference in the errors reported for $d=4$ and $d=5$ for these step sizes. Thus, for the rest of this work, we will use $d=4$ for all numerical computations of nonlocal Poisson equations.
\begin{table}[tb!]
	\centering
	\begin{tabularx}{\textwidth}{>{\centering\arraybackslash}c|>{\centering\arraybackslash}X|>{\centering\arraybackslash}X|>
			{\centering\arraybackslash}c|>
			{\centering\arraybackslash}X|>
			{\centering\arraybackslash}c|>
			{\centering\arraybackslash}X|>
			{\centering\arraybackslash}c}
		\hhline{|========|}
		&&&&&&&\\[-2ex]
		\hspace{7pt}$d$\hspace{7pt} & $h$ & $\overset{\beta=1.2}{\varepsilon_2}$ & Order & $\overset{\beta=2.0}{\varepsilon_2}$ & Order &  $\overset{\beta=2.5}{\varepsilon_2}$  & Order\\[0.5ex]
		\hline
		&&&&&&&\\[-2ex]
		4 & $2.00\cdot 10^{-2}$ & $1.29\cdot 10^{-04}$ & ----- &  $3.56\cdot 10^{-04} $ & ----- & $8.97\cdot 10^{-04}$ & -----\\
		& $1.00\cdot 10^{-2}$ & $2.00\cdot 10^{-06}$ & $6.0$ &  $4.10\cdot 10^{-06}$& $6.4$ & $1.09\cdot 10^{-05}$ & $6.4$ \\
		& $5.00\cdot 10^{-3}$ & $2.65\cdot 10^{-08}$ & $6.2$ &  $7.75\cdot 10^{-08}$ & $5.7$ & $2.32\cdot 10^{-07}$ & $5.6$\\
		& $2.50\cdot 10^{-3}$ & $1.44\cdot 10^{-10}$ & $7.5$ &  $2.19\cdot 10^{-10}$& $8.5$ & $1.61\cdot 10^{-09}$ & $7.2$ \\
		& $1.25\cdot 10^{-3}$ & $1.62\cdot 10^{-12}$ & $6.5$ &  $4.79\cdot 10^{-12}$ & $5.5$ & $2.24\cdot 10^{-11}$ & $6.2$\\[0.5ex]
		\hline
		&&&&&&&\\[-2ex]
		5 & $2.00\cdot 10^{-2}$ & $6.65\cdot 10^{-04}$ & ----- & $2.09\cdot 10^{-03}$ & ----- & $5.11\cdot 10^{-03} $ & -----\\
		& $1.00\cdot 10^{-2}$ & $6.66\cdot 10^{-06}$ & $6.6$ & $1.71\cdot 10^{-05}$ & $6.9$ & $4.99\cdot 10^{-05}$ & $6.7$\\
		& $5.00\cdot 10^{-3}$ & $4.94\cdot 10^{-08}$ & $7.1$ & $1.39\cdot 10^{-07}$ & $6.9$ & $4.17\cdot 10^{-07}$ & $6.9$\\
		& $2.50\cdot 10^{-3}$ & $2.21\cdot 10^{-10}$ & $7.8$ & $5.84\cdot 10^{-10}$ & $7.9$ & $2.63\cdot 10^{-09}$ & $7.3$\\
		& $1.25\cdot 10^{-3}$ & $1.31\cdot 10^{-12}$ & $7.4$ & $4.41\cdot 10^{-12}$ & $7.0$  & $2.06\cdot 10^{-11}$ & $7.0$\\[0.5ex]
		\hhline{|========|}
	\end{tabularx}
	\caption{Convergence of the 2D-FC based solution of the nonlocal Poisson equation \eqref{nonlocal_poisson} in the settings of Example~\ref{convergence_test_poisson} with fixed $\delta=0.4$, fixed $M=d+1$, and two different values of $d$.}
	\label{conv_study_nonlocal_Poisson}
\end{table}

\subsection{Accuracy and convergence test for nonlocal diffusion equation}
Similarly, we construct a manufactured solution to study the accuracy and convergence for the nonlocal diffusion equation~\eqref{nonlocal_diffusion_eqn}. Let
\begin{equation*}
	u(x,y,t) = \sin(2\pi r_1x)\sin(2\pi r_2y)\exp({-2\pi^2\kappa t})
\end{equation*}
be the solution of \eqref{nonlocal_diffusion_eqn}, where $\kappa=0.1$, $r_1=r_2=15.6455$. Clearly, $u(x,y)$ is an oscillatory non-periodic function. To obtain the heat source, $s(x,y,t)$, we apply operator $\Ldel$ on $u$, using the property of $\Ldel$ given in \eqref{multipliers_property},
with fixed $\delta=0.3$ and three values of $\beta: 1.0$, $2.0$, and $2.5$. Thus,
\begin{equation}
	s(x,y,t) = u_t-\Ldel u = \left(-2\pi^2\kappa - \mdel\right) u.\label{manufactured_heat_source}
\end{equation}
A direct computation of the multipliers $\mdel$ is provided in Table~\ref{manufactured_sol_multipliers_heat}. We compute these multipliers using the implementation of the $_2F_3$ hypergeometric function provided by mpmath \cite{johansson2013mpmath}, as detailed in \cite{alali2021fourier}.
\begin{table}[ht!]
	\centering
	\begin{tabularx}{\textwidth}{>{\centering\arraybackslash}c|>
			{\centering\arraybackslash}X|>
			{\centering\arraybackslash}X|>
			{\centering\arraybackslash}X}
		\hhline{|====|}
		&&&\\[-2ex]
		\hspace{7pt}$\delta$\hspace{7pt} & $\overset{\beta = 1.0}{\mdel}$ & $\overset{\beta=2.0}{\mdel}$ & $\overset{\beta=2.5}{\mdel}$ \\[0.5ex]
		\hline
		&&&\\[-2ex]
		$0.3$& $-130.16228859689554$ & $-321.3202730766787$ & $-689.8419741563309$\\
		[0.5ex]
		\hhline{|====|}
	\end{tabularx}
	\caption{The multipliers of the nonlocal operator $\Ldel$ applied on the function $u=\sin(2\pi r_1 x)\sin(2\pi r_2 y)$, where $r_1= r_2 = 15.6455$, with a fixed $\delta=0.3$ and three different values of $\beta$.}
	\label{manufactured_sol_multipliers_heat}
\end{table}

\begin{ex}[Convergence test for nonlocal diffusion equation]\label{convergence_test_diffusion}
	We employ our proposed method to study the error convergence of the nonlocal diffusion problem $\eqref{nonlocal_diffusion_eqn}$ over the kite-shaped domain contained within the curve given by \eqref{kite_domain}
	using step sizes $h=0.02/2^k,~k=0,1,2,3,4$ and time step size $\tau=5.0\cdot 10^{-07}$. We fix $\delta=0.3$, and choose three values of $\beta:$ $1.0$, $2.0$, and $2.5$, and we set the boundary function $b(x,y) =  \sin(2\pi r_1x)\sin(2\pi r_2y)\exp({-2\pi^2\kappa t})$, where $\kappa=0.1$, $r_1=r_2=15.6455$, and the heat source, $s(x,y,t)$, is as given in \eqref{manufactured_heat_source}.
\end{ex}

A convergence study for this test case is shown in Table~\ref{conv_study_nonlocal_diffusion} for $d=4$ and $d=5$, with fixed $M=d+1$. We use the manufactured solution described above to compute the $L^2$ errors reported. We observe sixth and seventh orders of convergence for $d=4$ and $d=5$, respectively, showing convergence of order $d+2$ for the three cases of $\beta$ considered. 
Moreover, as observed in the case of nonlocal Poisson equation, we observe no significant difference in the errors reported for $d=4$ and $d=5$. Thus, for the rest of this work, we will use $d=4$ for all numerical computations of nonlocal diffusion equations.

\begin{table}[tb!]
	\centering
	\begin{tabularx}{\textwidth}{>{\centering\arraybackslash}c|>{\centering\arraybackslash}X|>{\centering\arraybackslash}X|>
			{\centering\arraybackslash}c|>
			{\centering\arraybackslash}X|>
			{\centering\arraybackslash}c|>
			{\centering\arraybackslash}X|>
			{\centering\arraybackslash}c}
		\hhline{|========|}
		&&&&&&&\\[-2ex]
		\hspace{7pt}$d$\hspace{7pt} & $h$ & $\overset{\beta=1.0}{\varepsilon_2}$ & Order & $\overset{\beta=2.0}{\varepsilon_2}$ & Order &  $\overset{\beta=2.5}{\varepsilon_2}$  & Order\\[0.5ex]
		\hline
		&&&&&&&\\[-2ex]
		4 & $2.00\cdot 10^{-2}$ & $1.13\cdot 10^{-05}$ & ----- & $5.33\cdot 10^{-05}$ & ----- & $2.48\cdot 10^{-04}$ & -----\\
		& $1.00\cdot 10^{-2}$ & $1.51\cdot 10^{-07}$ & $6.2$ & $9.45\cdot 10^{-07}$ & $5.8$ & $6.81\cdot 10^{-06}$ & $5.2$\\
		& $5.00\cdot 10^{-3}$ & $2.07\cdot 10^{-09}$ & $6.2$ & $1.48\cdot 10^{-08}$ & $6.0$ & $1.45\cdot 10^{-07}$ & $5.6$ \\
		& $2.50\cdot 10^{-3}$ & $1.87\cdot 10^{-11}$ & $6.8$ & $9.21\cdot 10^{-11}$ & $7.3$ & $1.22\cdot 10^{-09}$ & $6.9$ \\
		& $1.25\cdot 10^{-3}$ & $2.43\cdot 10^{-13}$ & $6.3$ & $1.29\cdot 10^{-12}$ & $6.2$ & $2.29\cdot 10^{-11}$ & $5.7$ \\[0.5ex]
		\hline
		&&&&&&&\\[-2ex]
		5 & $2.00\cdot 10^{-2}$ & $7.87\cdot 10^{-05}$ & ----- & $4.10\cdot 10^{-04}$ & ----- & $1.93\cdot 10^{-03}$ & -----\\
		& $1.00\cdot 10^{-2}$ & $8.84\cdot 10^{-07}$ & $6.5$ & $5.86\cdot 10^{-06}$ & $6.1$ & $4.16\cdot 10^{-05}$ & $5.5$\\
		& $5.00\cdot 10^{-3}$ & $5.53\cdot 10^{-09}$ & $7.3$ & $3.30\cdot 10^{-08}$ & $7.5$ & $3.23\cdot 10^{-07}$ & $7.0$\\
		& $2.50\cdot 10^{-3}$ & $2.99\cdot 10^{-11}$ & $7.5$ & $1.87\cdot 10^{-10}$ & $7.5$ & $2.49\cdot 10^{-09}$ & $7.0$\\
		& $1.25\cdot 10^{-3}$ & $3.49\cdot 10^{-13}$ & $6.4$ & $1.27\cdot 10^{-12}$ & $7.2$ & $2.15\cdot 10^{-11}$ & $6.9$ \\[0.5ex]
		\hhline{|========|}
	\end{tabularx}
	\caption{Convergence of the 2D-FC based solution of the nonlocal diffusion equation \eqref{nonlocal_diffusion_eqn} in the settings of Example~\ref{convergence_test_diffusion} with fixed $\delta=0.3$, fixed $M=d+1$, two different values of $d$: $4$ and $5$, and three choices of $\beta$: $1.0$, $2.0$, and $2.5$.}
	\label{conv_study_nonlocal_diffusion}
\end{table}

\subsection{Study of discontinuities for the nonlocal Poisson equation}
In this section, we study the behavior of the discontinuities in the numerical solutions of nonlocal Poisson equations on bounded domains. For the case of integrable kernels ($\beta<2$), we emphasize that if $f$ is discontinuous, then the solution of the nonlocal Poisson equation is also discontinuous and will share the same location of jumps with $f$. However, for the case of singular kernels with $\beta>3$, we form a conjecture that even if $f$ is discontinuous, the solution will be continuous.

\begin{ex}\label{no_discont}
	Let $\Omega$ be the kite-shaped domain contained within the curve defined in \eqref{kite_domain}. We solve the nonlocal Poisson equation
	\begin{eqnarray}
		\begin{cases}
			\Ldel u(x,y) = 4 ,& (x,y)\in \Omega \\
			u(x,y) = x^2+y^2, & (x,y)\in \collar
		\end{cases}
	\end{eqnarray}
	with $\delta=0.4$, three different values of $\beta:$ $1.0$, $2.0$, $3.0$, and we choose the 2D-FC parameters $d=4$, $M=5$, and $h=k_1=k_2=0.0025$.
\end{ex}

\begin{figure}[ht!]
    \centering
    \begin{subfigure}[h]{0.28\textwidth}
        \centering
        \includegraphics[width=\textwidth]{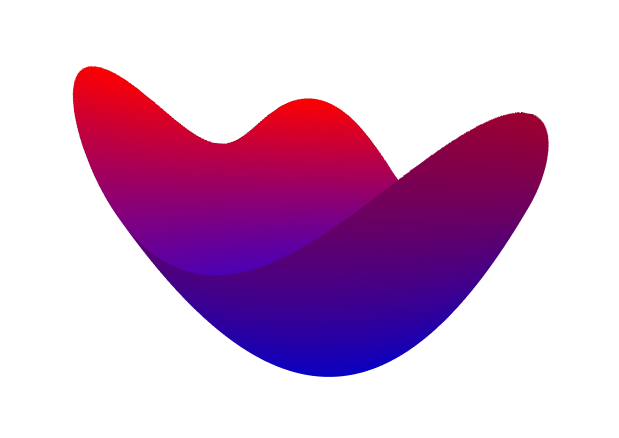}
        \caption{$\beta=1.0$}
    \end{subfigure}
    \hfill
    \begin{subfigure}[h]{0.28\textwidth}
        \centering
        \includegraphics[width=\textwidth]{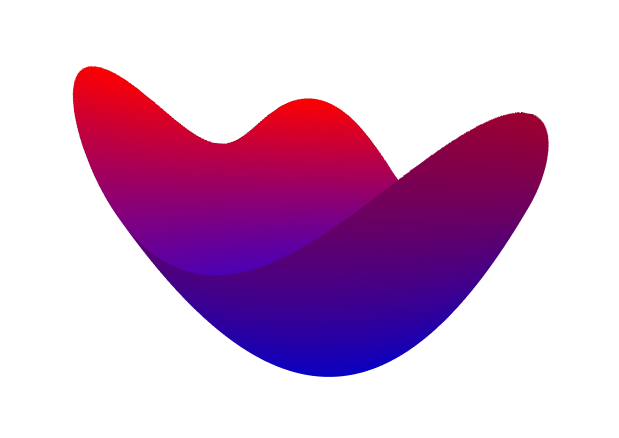}
        \caption{$\beta=2.0$}
    \end{subfigure}
    \hfill
    \begin{subfigure}[h]{0.28\textwidth}
        \centering
        \includegraphics[width=\textwidth]{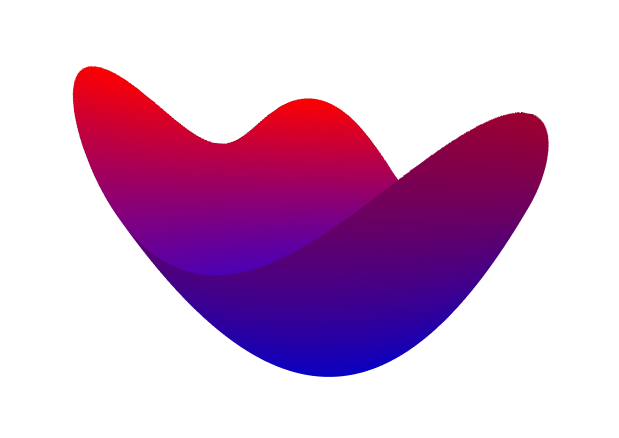}
        \caption{$\beta=3.0$}
    \end{subfigure}
    \caption{Simulation results in Example~\ref{no_discont} with a fixed $\delta=0.4$ and three different values of $\beta$.}
    \label{Poisson_sim2}
\end{figure}

In this example, since for any $\delta>0$ and $\beta<4$, $\Ldel(x^2+y^2)=4$ for all $x\in \Omega$, then the unique solution for this problem is given explicitly by $u(x,y)=x^2+y^2$, which is obviously a continuous function.  
The numerical solution obtained using our solver is consistent with the above observation as illustrated in Figure~\ref{Poisson_sim2}.

\begin{figure}[t!]
    \centering
    \begin{subfigure}[h]{0.24\textwidth}
        \centering
        \includegraphics[width=\textwidth]{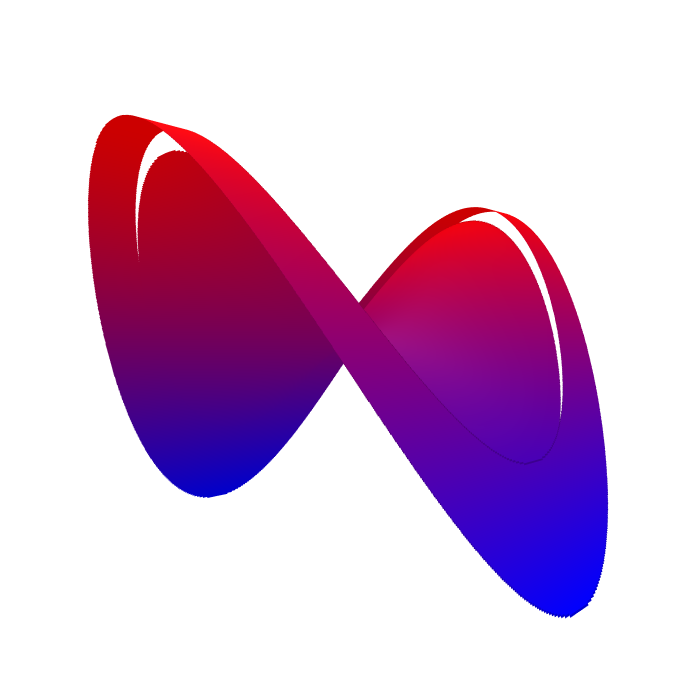}
        \caption{$\beta=1.5$}
    \end{subfigure}
    \hfill
    \begin{subfigure}[h]{0.24\textwidth}
        \centering
        \includegraphics[width=\textwidth]{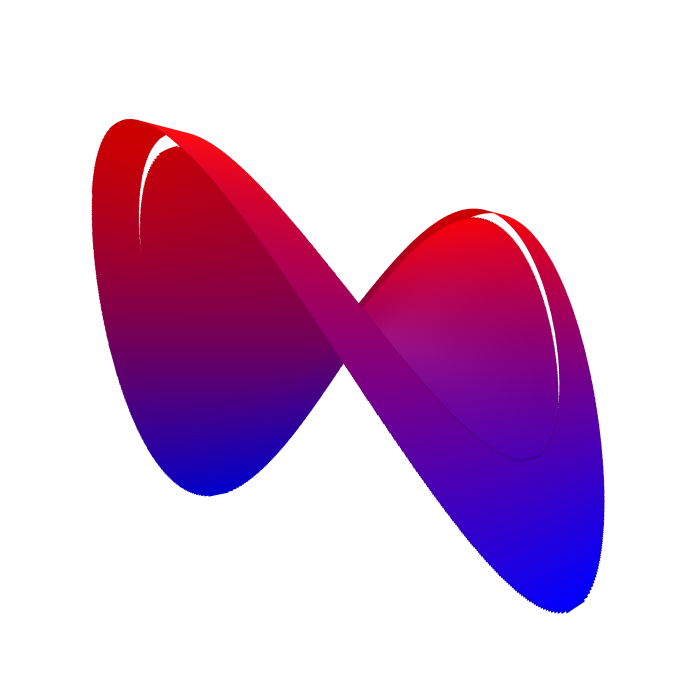}
        \caption{$\beta=2.0$}
    \end{subfigure}
    \hfill
    \begin{subfigure}[h]{0.24\textwidth}
        \centering
        \includegraphics[width=\textwidth]{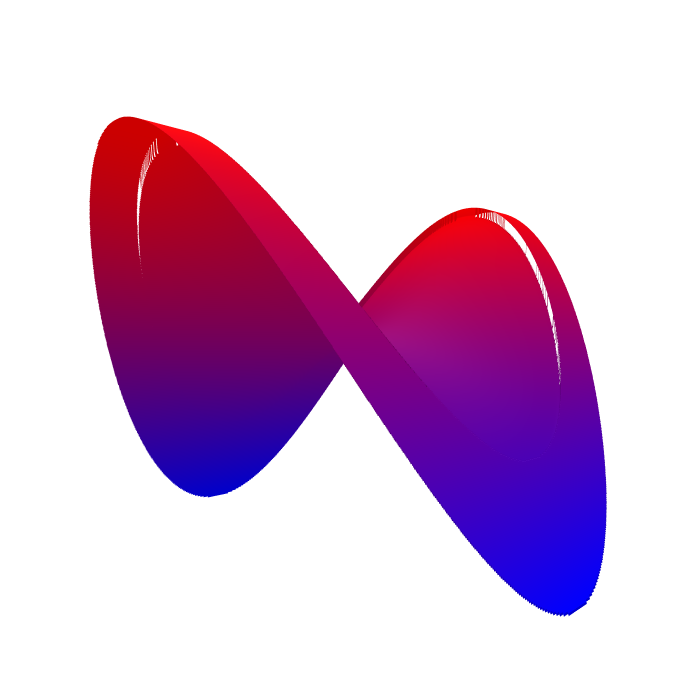}
        \caption{$\beta=2.2$}
    \end{subfigure}
    \hfill
    \begin{subfigure}[h]{0.24\textwidth}
        \centering
        \includegraphics[width=\textwidth]{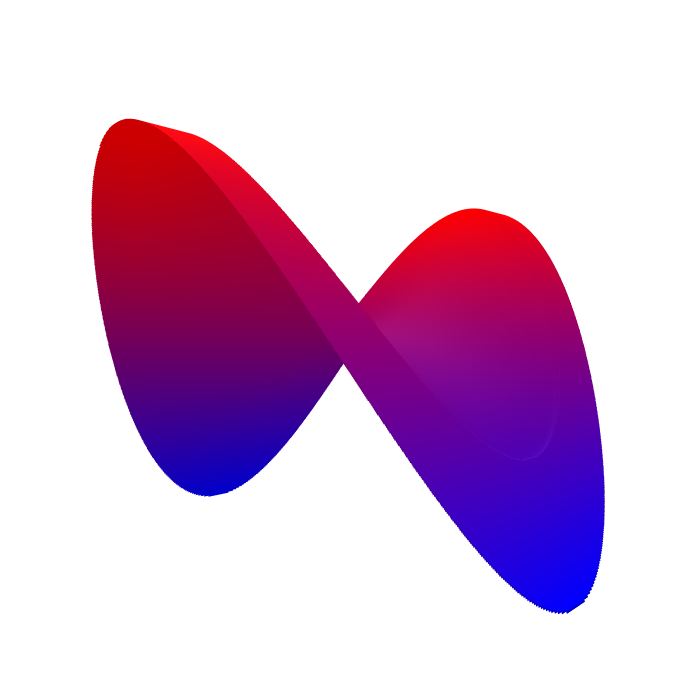}
        \caption{$\beta=3.1$}
    \end{subfigure}
    \caption{Simulation results in Example~\ref{max_bdy_val} with a fixed $\delta=0.2$ and four different values of $\beta$. We observe that the magnitude of the jump at the $\partial \Omega$ decreases as $\beta$ increases (see Table~\ref{magnitude_of_jump1}).}
    \label{Poisson_sim1}
\end{figure}

\begin{ex}\label{max_bdy_val}
	Let $\Omega$ be the unit disk. We solve the nonlocal Poisson equation 
	\begin{eqnarray}
		\begin{cases}
			\Ldel u(x,y) = |\sin(\pi x)\sin(\pi y)|,& (x,y)\in \Omega \\
			u(x,y) = -x^2, & (x,y)\in \collar
		\end{cases}
	\end{eqnarray}
	with $\delta=0.2$, four values of $\beta:$ $1.5$, $2.0$, $2.2$, $3.1$, and we choose the 2D-FC parameters $d=4$, $M=5$, and $h=k_1=k_2=0.0025$.
\end{ex}
Figure~\ref{Poisson_sim1} presents the numerical solution for the three choices of $\beta$ considered. We observe that the solution is discontinuous at $\partial \Omega$, and the magnitude of the jump decreases as $\beta$ increases. Specifically, we observe that for $\beta>3$, the magnitude of the jump is small. The maximum and the minimum magnitudes of the jumps at $\partial \Omega$ for all the values of $\beta$ considered are provided in Table~\ref{magnitude_of_jump1}, and they are computed using the following formula:
\begin{equation*}
	\text{magnitude of the jump} = \left|\widetilde{u}(x_\theta, y_\theta) - b(x_\theta, y_\theta)\right|,
\end{equation*}
where $\widetilde{u}$ is the extension of $u$ from the interior of $\Omega$ to the boundary $\partial \Omega$, $x_\theta$ and $y_\theta$ are the $x$ and $y$ coordinates of the boundary $\partial \Omega$ at $\theta\in\Theta_B$, and $b$ is the given boundary function.

\begin{table}[t!]
    \centering
        \begin{tabularx}{\textwidth}{>{\centering\arraybackslash}c|>
        {\centering\arraybackslash}X|>
        {\centering\arraybackslash}X|>
        {\centering\arraybackslash}X|>
        {\centering\arraybackslash}X}
        \hhline{|=====|}
        &&&\\[-2ex]
        \hspace{7pt}\hspace{7pt} & $\beta = 1.5$ & $\beta=2.0$ & $\beta=2.2$ & $\beta=3.1$ \\[0.5ex]
        \hline
        &&&\\[-2ex]
        $\min$ magnitude& $0.0538468$ & $0.0286598$ & $0.0187264$ & $6.24\cdot 10^{-05}$\\
        $\max$ magnitude& $0.0646947$ & $0.0454215$ & $0.0372901$ & $0.0097378$\\[0.5ex]
       \hhline{|=====|}
    \end{tabularx}
    \caption{The minimum and the maximum magnitudes of the jumps at $\partial \Omega$ in the simulation results presented in Figure~\ref{Poisson_sim1}. We observe that the magnitudes of the jumps decrease as $\beta$ increases.}
    \label{magnitude_of_jump1}
\end{table}

The following example considers the case where $f(x,y)$ is discontinuous. As in the periodic case, we expect the solution to be discontinuous at the location where $f$ is discontinuous.

\begin{ex}\label{discont_rhs}
	Let $\Omega$ be the kite-shaped domain contained within the curve defined in \eqref{kite_domain} and let
	\[f(x,y) = \begin{cases}
		80,& \text{ if } x^2+4y^2<0.2\\
		4, & \text{ otherwise}
	\end{cases}.\] 
	We solve the nonlocal Poisson equation 
	\begin{eqnarray}
		\begin{cases}
			\Ldel u(x,y) = f(x,y),& (x,y)\in \Omega \\
			u(x,y) = x^2+y^2, & (x,y)\in \collar
		\end{cases}
	\end{eqnarray}
	with $\delta=0.5$, four different values of $\beta:$ $1.2$, $2.0$, $2.5$, $3.1$, and we choose the 2D-FC parameters $d=4$, $M=5$, and $h=k_1=k_2=0.0025$. 
\end{ex}

\begin{figure}[b!]
    \centering
    \begin{subfigure}[h]{0.24\textwidth}
        \centering
        \includegraphics[width=\textwidth]{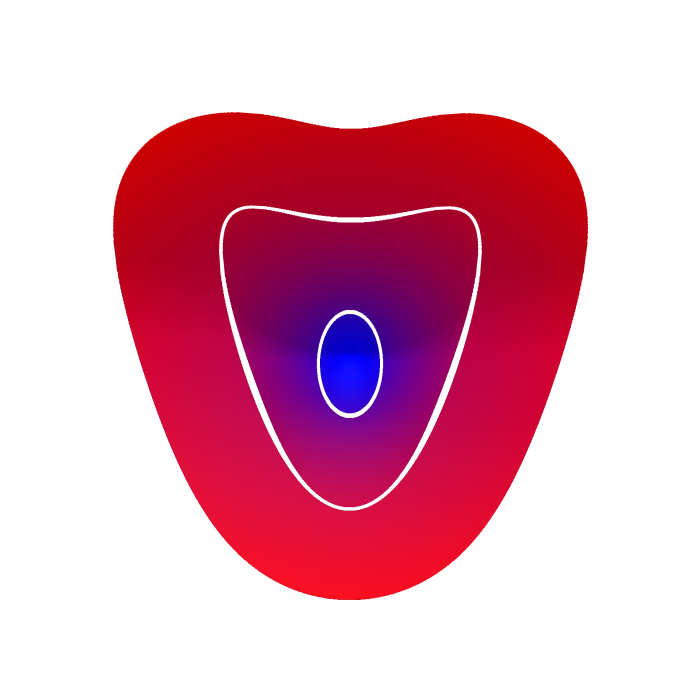}
        \caption{$\beta=1.2$}
    \end{subfigure}
    \hfill
    \begin{subfigure}[h]{0.24\textwidth}
        \centering
        \includegraphics[width=\textwidth]{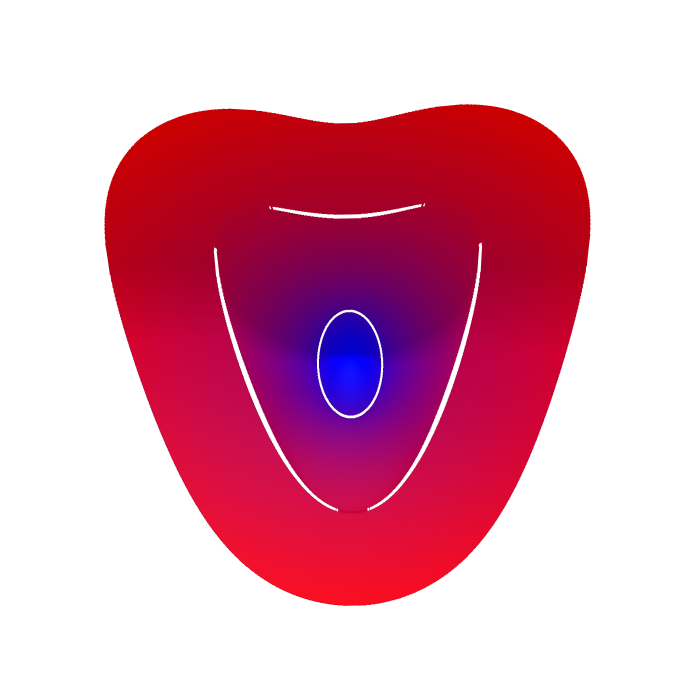}
        \caption{$\beta=2.0$}
    \end{subfigure}
    \hfill
    \begin{subfigure}[h]{0.24\textwidth}
        \centering
        \includegraphics[width=\textwidth]{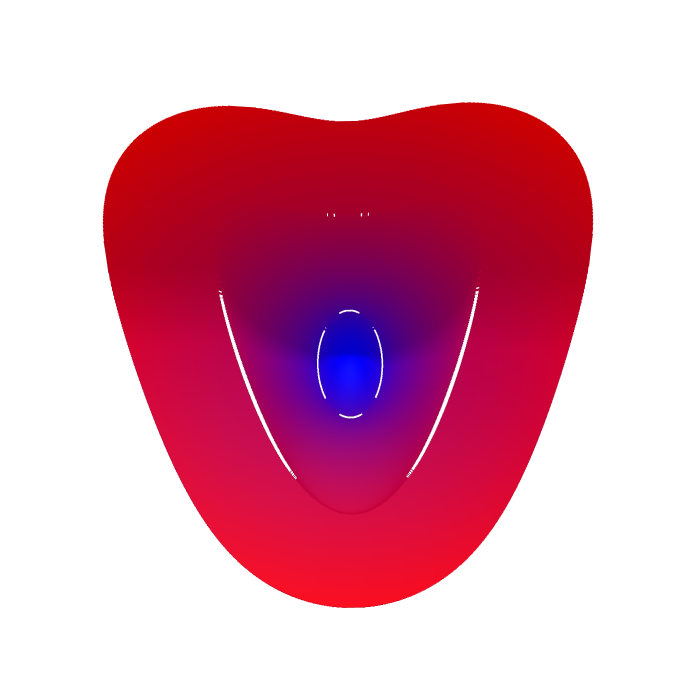}
        \caption{$\beta=2.5$}
    \end{subfigure}
    \hfill
    \begin{subfigure}[h]{0.24\textwidth}
        \centering
        \includegraphics[width=\textwidth]{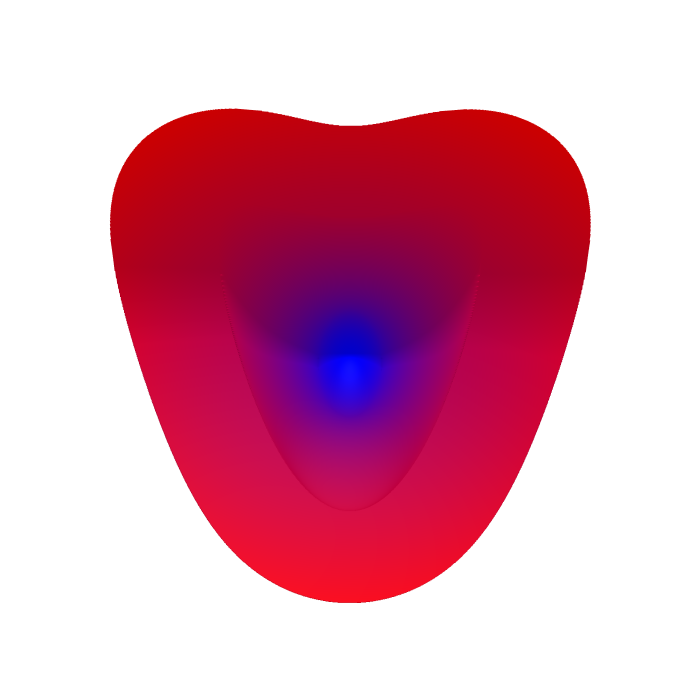}
        \caption{$\beta=3.1$}
    \end{subfigure}
    \caption{Simulation results in Example~\ref{discont_rhs} with a fixed $\delta=0.5$ and four different values of $\beta$.}
    \label{Poisson_sim3}
\end{figure}

In Figure~\ref{Poisson_sim3}, we present the numerical solutions for four different values of $\beta$: $1.2$, $2.0$, $2.5$, and $3.1$. We observe that as $\beta$ increases, the magnitudes of the jump discontinuity inside the domain $\Omega$ and at its boundary $\partial \Omega$ decrease. Furthermore, these discontinuities occur at the exact location where $f$ is discontinuous. Specifically, the maximum and the minimum magnitudes of the jumps in the interior of the domain along the ellipse for all the values of $\beta$ considered are given in Table~\ref{magnitude_of_jump2}, and they are computed using the following formula:
\begin{equation*}
	\text{magnitude of the jump} = \left|\widetilde{u}_{\text{int}}(x_\theta, y_\theta) - \widetilde{u}_{\text{ext}}(x_\theta, y_\theta)\right|,
\end{equation*}
where $\widetilde{u}_{\text{int}}$ is the extension of $u$ from the interior of the ellipse to the boundary of the ellipse,  $\widetilde{u}_{\text{ext}}$ is the extension of $u$ from the exterior of the ellipse to the boundary of the ellipse, and $x_\theta$ and $y_\theta$ are the $x$ and $y$ coordinates of the boundary of the ellipse at $\theta\in\Theta_B$.

\begin{table}[b!]
    \centering
        \begin{tabularx}{\textwidth}{>{\centering\arraybackslash}c|>
        {\centering\arraybackslash}X|>
        {\centering\arraybackslash}X|>
        {\centering\arraybackslash}X|>
        {\centering\arraybackslash}X}
        \hhline{|=====|}
        &&&&\\[-2ex]
        \hspace{7pt}\hspace{7pt} & $\beta = 1.2$ & $\beta=2.0$ & $\beta=2.5$ & $\beta=3.1$ \\[0.5ex]
        \hline
        &&&&\\[-2ex]
        $\min$ magnitude& $1.3589861$ & $0.3237825$ & $0.0363031$ & $1.75\cdot 10^{-05}$\\
        $\max$ magnitude& $1.3681906$ & $0.4040726$ & $0.0856388$ & $0.0067363$\\[0.5ex]
       \hhline{|=====|}
    \end{tabularx}
    \caption{The minimum and the maximum magnitudes of the jumps in the interior of $\Omega$ in the simulation results presented in Figure~\ref{Poisson_sim3}. We observe that the magnitudes of the jumps decrease as $\beta$ increases.}
    \label{magnitude_of_jump2}
\end{table}


To make a connection between the numerical results in this work, in particular Example \ref{discont_rhs}, and previous mathematical analysis in this area, we note that regularity results for the nonlocal Poisson equation on periodic domains have been developed in \cite{alali2021fourier}. However, regularity of solutions for the nonlocal equation on bounded domains is still an open problem. 

For the periodic case, Theorem~4.1 in \cite{alali2021fourier} implies that if $f\in L^2$, then $u\in H^{\beta-n}$. In particular, in the case when when $n=2$ and $3<\beta<4$, then $u\in H^{\beta-2}$, which implies that $u$ is continuous using the embedding theorem. Using intuition from the periodic case combined with the numerical observations in Example~\ref{discont_rhs}, we expect that this smoothing behavior as $\beta$ increases would be true in the case of bounded domains  as stated in the following conjecture. 
\begin{conjecture}
	Let $u$ be the solution of \eqref{nonlocal_poisson} and assume that $f\in L^2(\Omega)$. Then, if $3<\beta<4$, then $u\in C(\Omega)$. \label{conjecture:continuity}
\end{conjecture}

\subsection{Evolution of discontinuities for nonlocal diffusion equation}
In this section, we study the propagation of discontinuity in the solution of nonlocal diffusion equation~\eqref{nonlocal_diffusion_eqn}. 

\begin{figure}[ht!]
    \centering
    \begin{subfigure}[h]{0.24\textwidth}
        \centering
        \includegraphics[width=\textwidth]{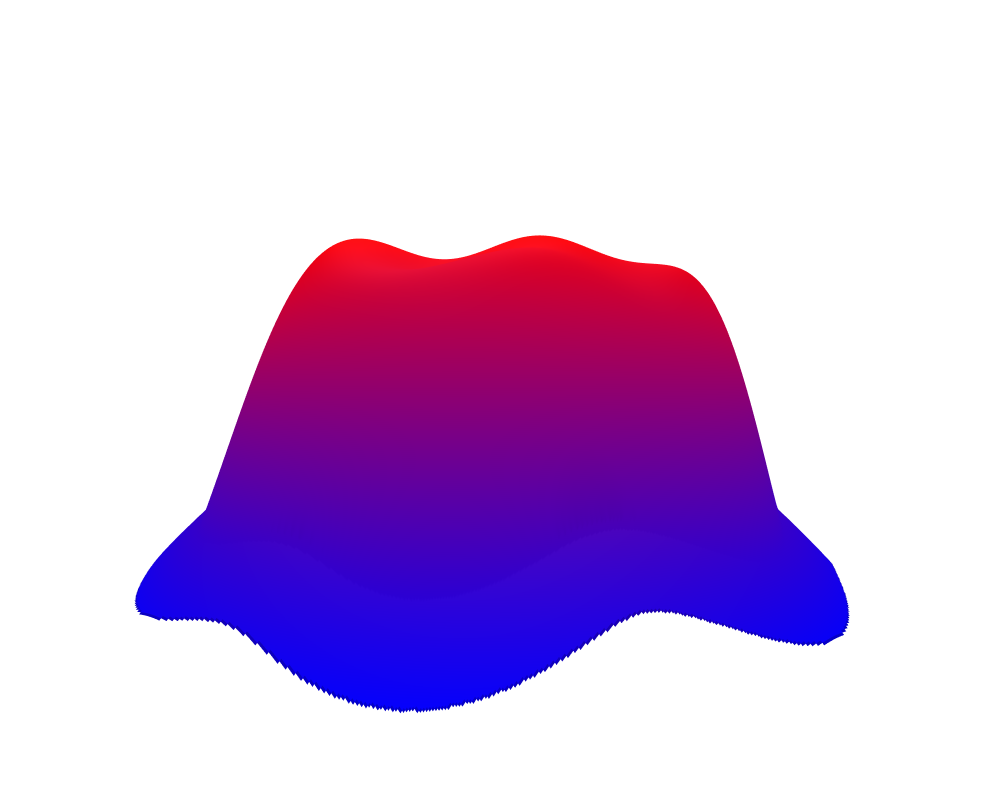}
    \end{subfigure}
    \hfill
    \begin{subfigure}[h]{0.24\textwidth}
        \centering
        \includegraphics[width=\textwidth]{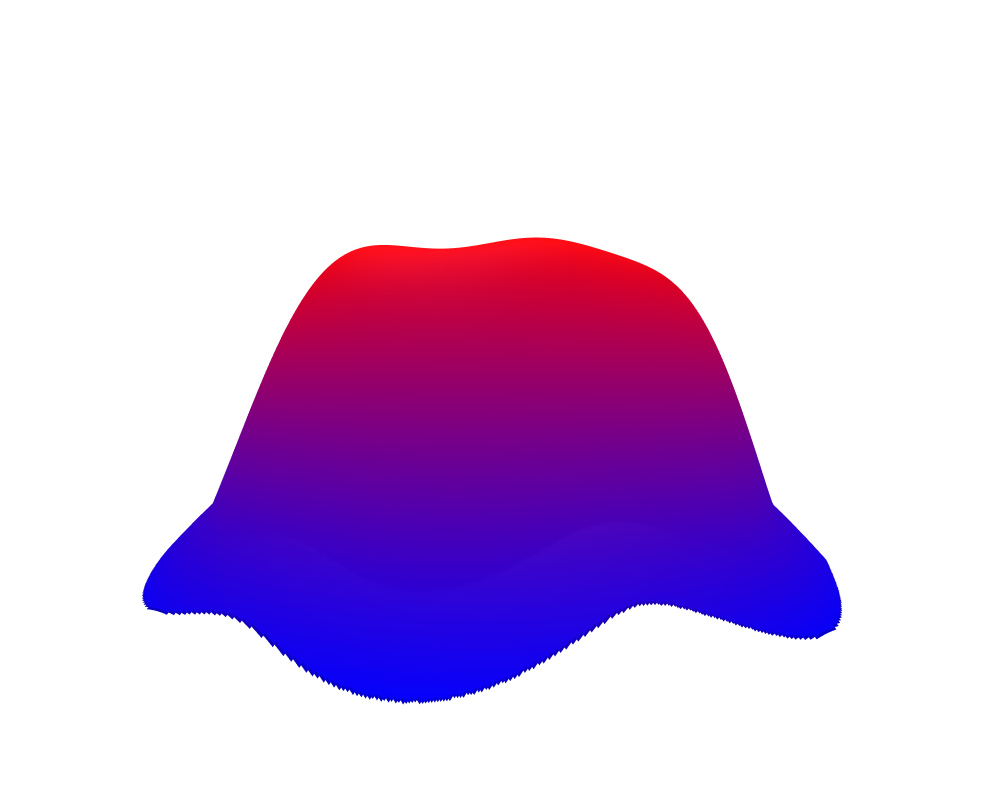}
    \end{subfigure}
    \hfill
    \begin{subfigure}[h]{0.24\textwidth}
        \centering
        \includegraphics[width=\textwidth]{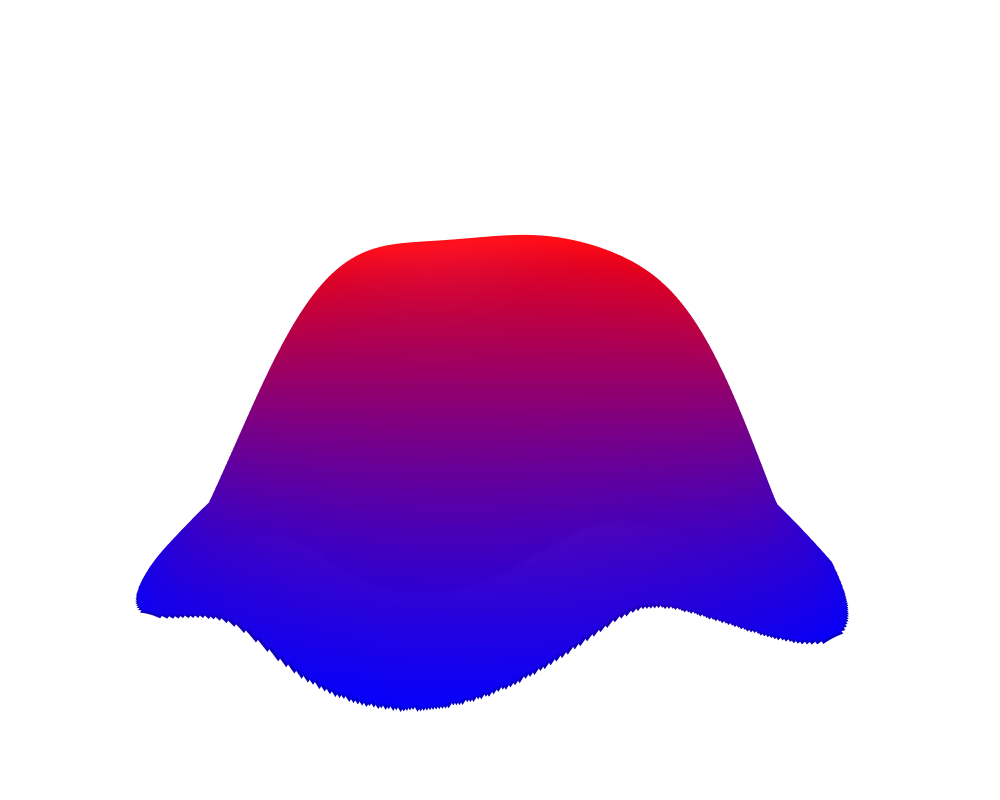}
    \end{subfigure}
    \hfill
    \begin{subfigure}[h]{0.24\textwidth}
        \centering
        \includegraphics[width=\textwidth]{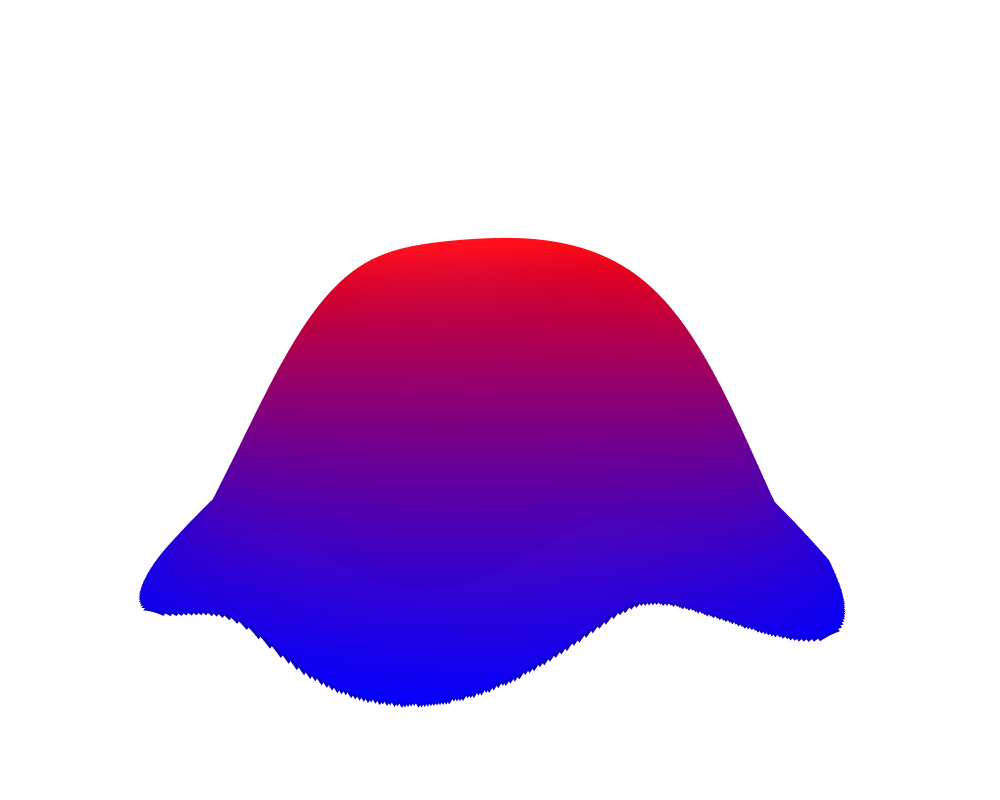}
    \end{subfigure}
    \begin{subfigure}[h]{0.23\textwidth}
        \centering
        \includegraphics[width=\textwidth]{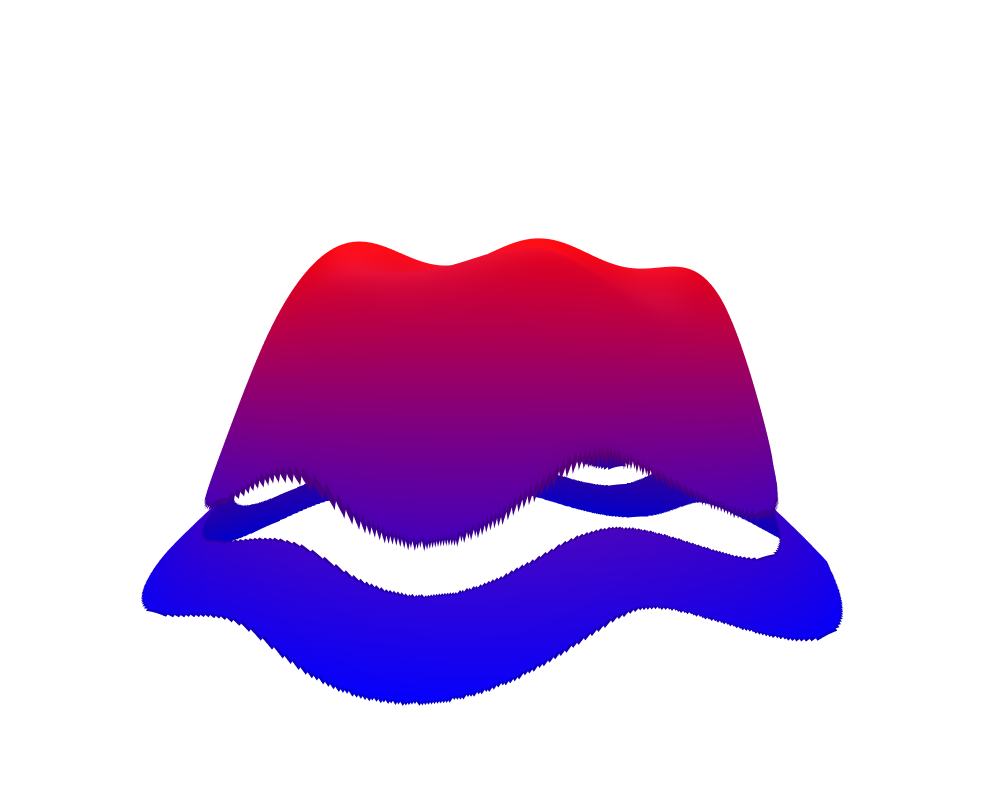}
        \caption{$t=0.02$}
    \end{subfigure}
    \hfill
    \begin{subfigure}[h]{0.24\textwidth}
        \centering
        \includegraphics[width=\textwidth]{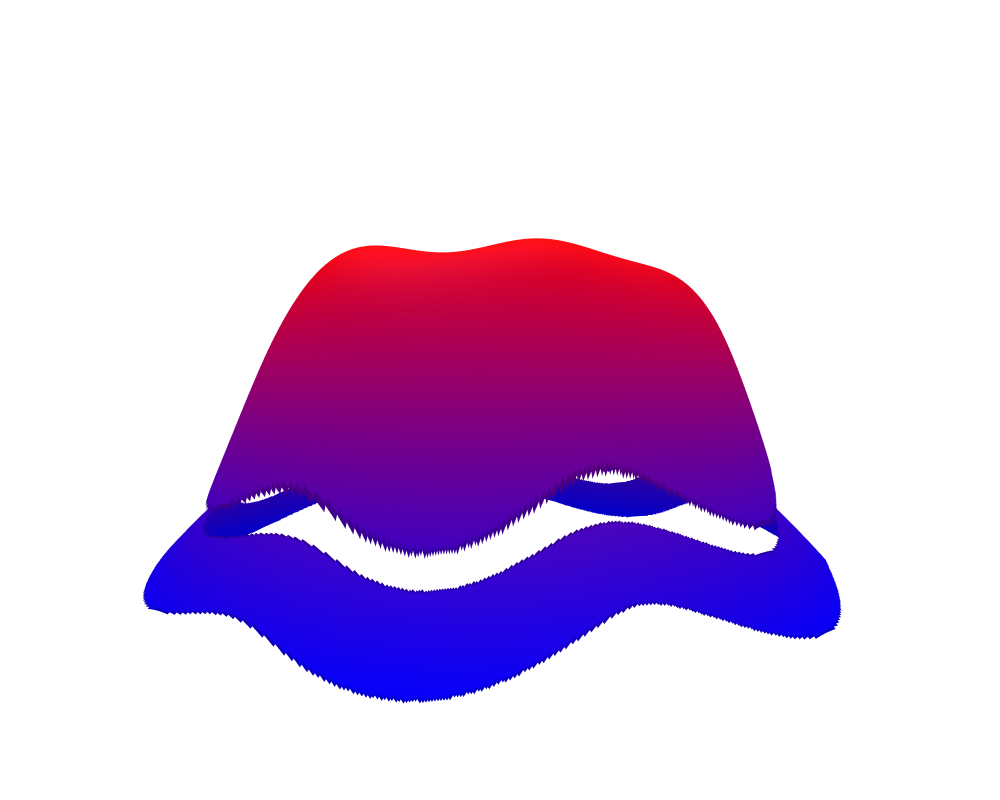}
        \caption{$t=0.03$}
    \end{subfigure}
    \hfill
    \begin{subfigure}[h]{0.24\textwidth}
        \centering
        \includegraphics[width=\textwidth]{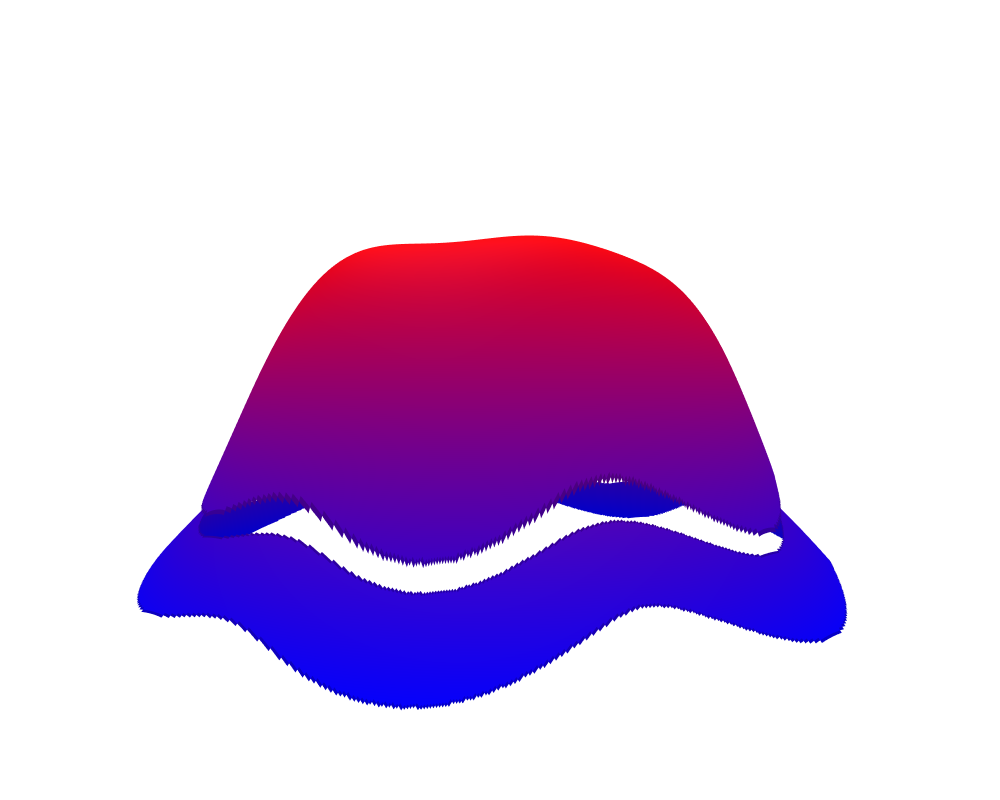}
        \caption{$t=0.04$}
    \end{subfigure}
    \hfill
    \begin{subfigure}[h]{0.24\textwidth}
        \centering
        \includegraphics[width=\textwidth]{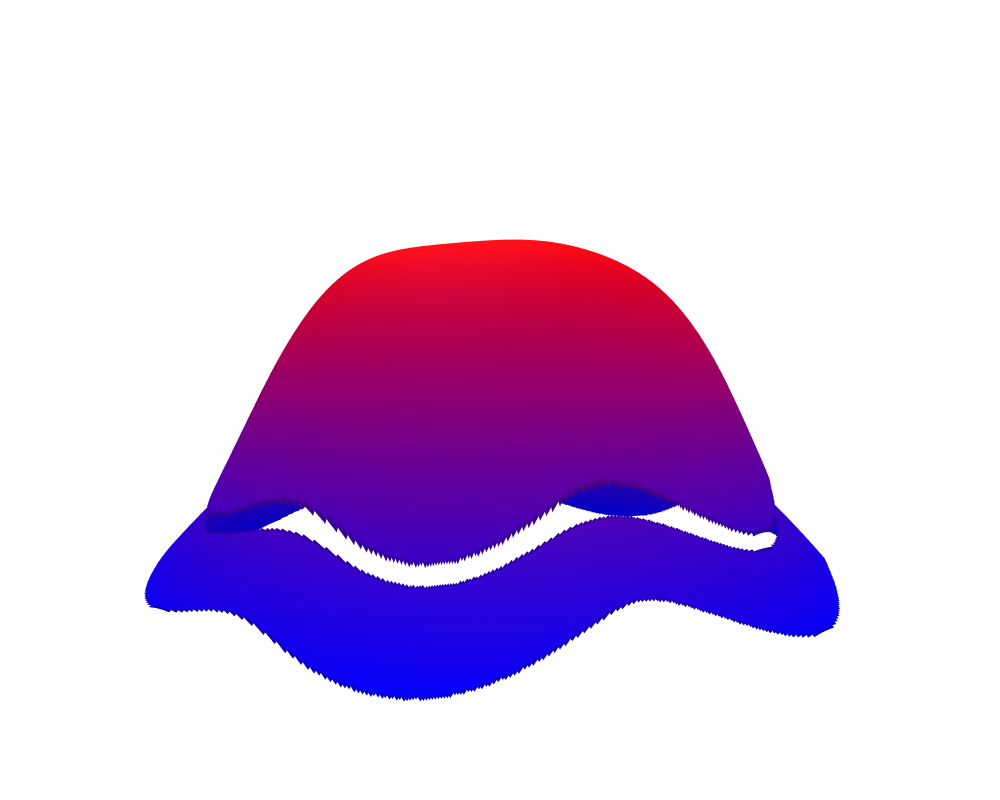}
        \caption{$t=0.05$}
    \end{subfigure}
    \caption{Simulation results in Example~\ref{nonlocal_diffusion1} with fixed $\delta=0.2$ and two different values of $\beta$: 4.0 (first row) and 1.5 (second row) at four different time instances.}
    \label{heat_sim1}
\end{figure}

\begin{ex}\label{nonlocal_diffusion1}
	Let $\Omega$ be the domain contained within the curve given by 
	\begin{eqnarray}
		x(\theta) = r(\theta)\cos(\theta); ~y(\theta) = r(\theta)\sin(\theta),
		\label{ring_shape_domain}
	\end{eqnarray}
	where $r(\theta)=1.1 + \cos(7\theta)/20 + \sin(4\theta)/30,$ and $0\le \theta\le 2\pi$.
	We solve the nonlocal diffusion equation
	\begin{eqnarray}
		\begin{cases}
			u_t(x,y,t) = \Ldel u(x,y,t),& (x,y,t)\in \Omega\times (0,T) \\
			u(x,y,0) = \sin(2\pi x)\cos(2\pi y),&(x,y)\in \Omega\\
			u(x,y,t) = -(1+x^2+y^2), & (x,y)\in \collar
		\end{cases}
	\end{eqnarray}
	with $\delta=0.2$ and two different values of $\beta$: $4.0$ and $1.5$ at four different time instances: $0.02$, $0.03$, $0.04$, and $0.05$ with time step size $\tau=2.5\cdot 10^{-7}$. We choose the 2D-FC parameters $d=4$, $M=5$, and $h=k_1=k_2=0.005$. 
\end{ex}
In this case, the corresponding  classical heat equation is given by 
\begin{eqnarray}\label{eq:classical_diffusion}
	\begin{cases}
		u_t(x,y,t) = \Delta u(x,y,t),& (x,y,t)\in \Omega\times (0,T), \\
		u(x,y,0) = \sin(2\pi x)\cos(2\pi y),&(x,y)\in \Omega,\\
		u(x,y) = -(1+x^2+y^2), & (x,y)\in \partial\Omega.
	\end{cases}
\end{eqnarray}
The first row of Figure~\ref{heat_sim1} demonstrates the numerical solution to the classical diffusion problem \eqref{eq:classical_diffusion} at different times $t=0.02$, $0.03$, $0.04$, and $0.05$.
The  second row of Figure~\ref{heat_sim1},  demonstrates the numerical solution for the nonlocal diffusion problem in Example \ref{nonlocal_diffusion1} with fixed $\delta=0.2$ and  $\beta=1.5$  for the same corresponding  times $t=0.02$, $0.03$, $0.04$, and $0.05$. For classical diffusion, we observe an instantaneous smoothing effect at $\partial\Omega$, in particular, the solution is continuous at the boundary.  In contrast,  the solution to the nonlocal problem is discontinuous at $\partial\Omega$, however, the magnitude of the jump decreases as time increases. 
We further observe that for both local and nonlocal diffusion, oscillations in the initial data smooth out  as time increases.

\begin{figure}[ht!]
    \centering
    \begin{subfigure}[h]{0.22\textwidth}
        \centering
        \includegraphics[width=\textwidth]{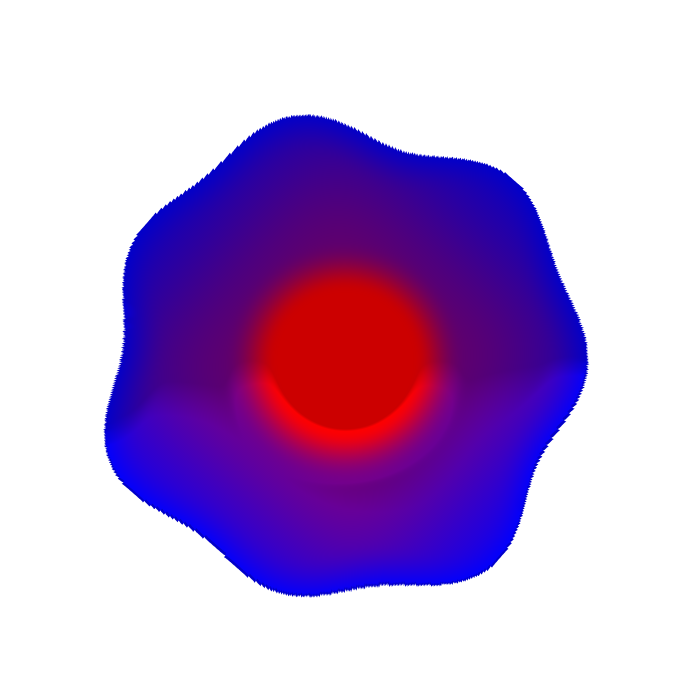}
    \end{subfigure}
    \hfill
    \begin{subfigure}[h]{0.22\textwidth}
        \centering
        \includegraphics[width=\textwidth]{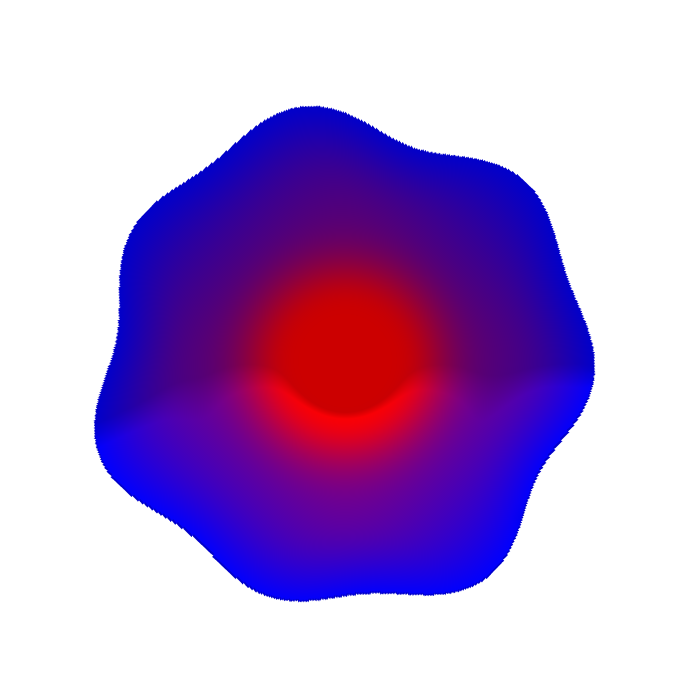}
    \end{subfigure}
    \hfill
    \begin{subfigure}[h]{0.22\textwidth}
        \centering
        \includegraphics[width=\textwidth]{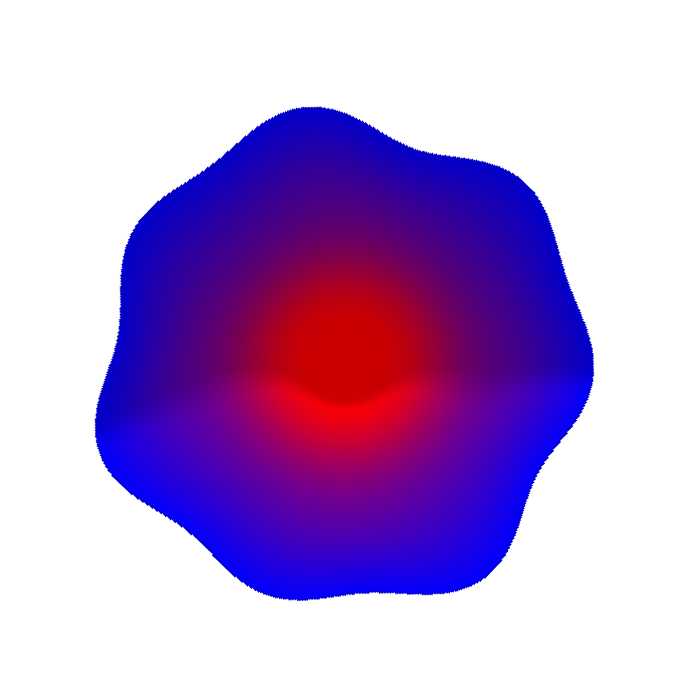}
    \end{subfigure}
    \hfill
    \begin{subfigure}[h]{0.22\textwidth}
        \centering
        \includegraphics[width=\textwidth]{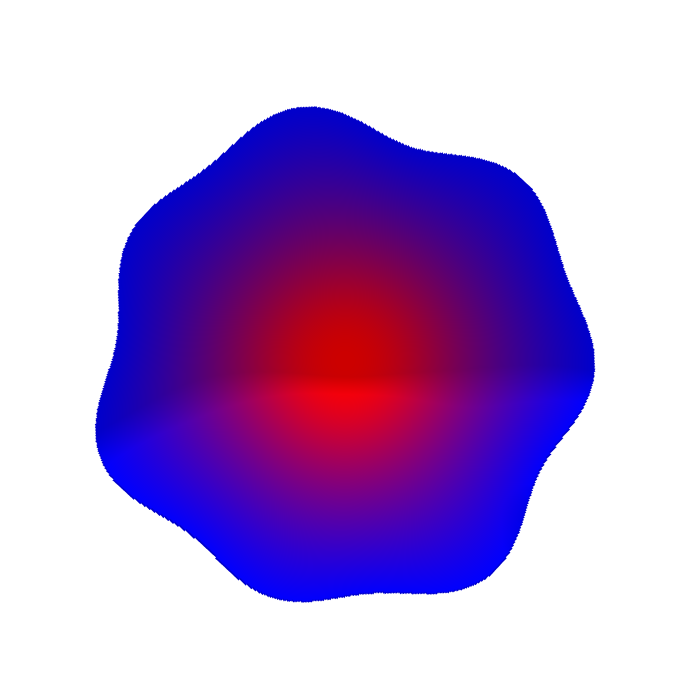}
    \end{subfigure}
    \begin{subfigure}[h]{0.22\textwidth}
        \centering
        \includegraphics[width=\textwidth]{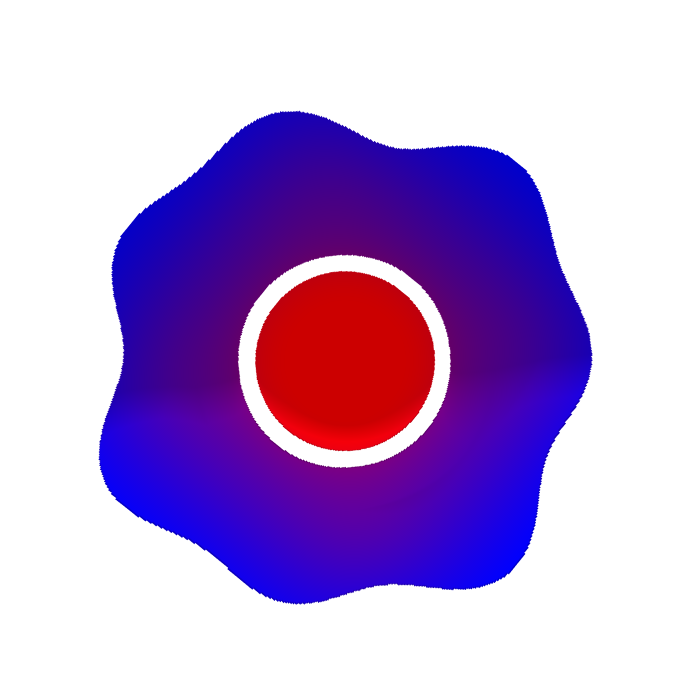}
        \caption{$t=0.001$}
    \end{subfigure}
    \hfill
    \begin{subfigure}[h]{0.22\textwidth}
        \centering
        \includegraphics[width=\textwidth]{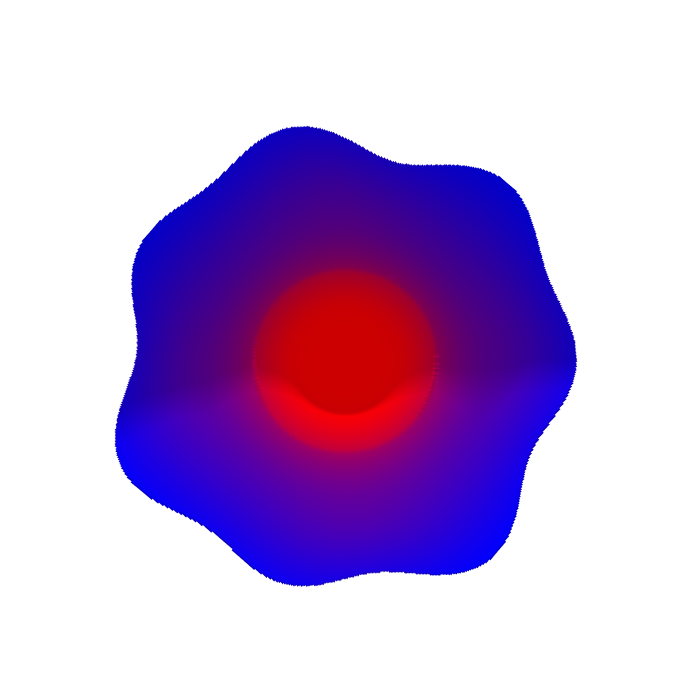}
        \caption{$t=0.005$}
    \end{subfigure}
    \hfill
    \begin{subfigure}[h]{0.22\textwidth}
        \centering
        \includegraphics[width=\textwidth]{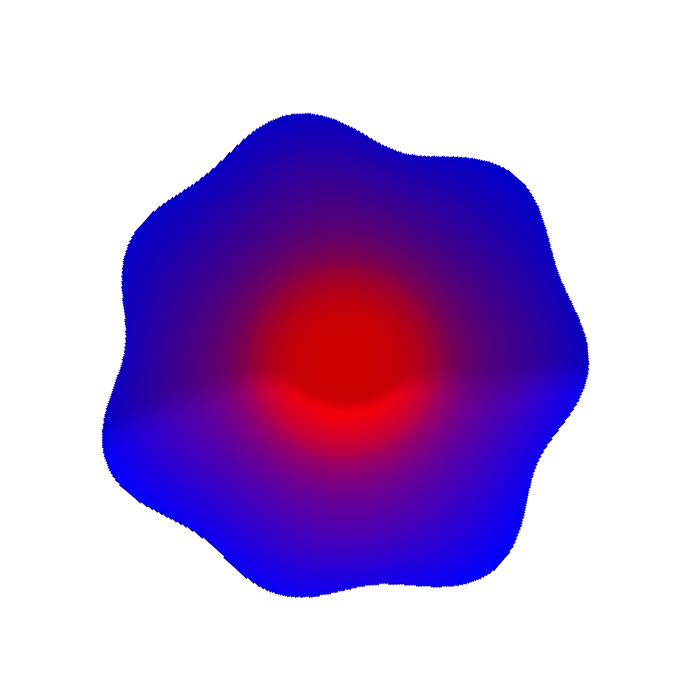}
        \caption{$t=0.01$}
    \end{subfigure}
    \hfill
    \begin{subfigure}[h]{0.22\textwidth}
        \centering
        \includegraphics[width=\textwidth]{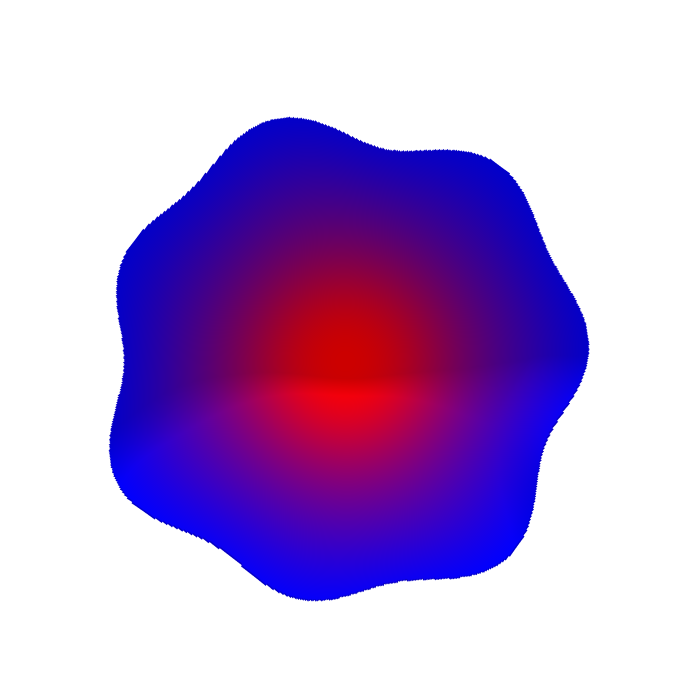}
        \caption{$t=0.03$}
    \end{subfigure}
    \caption{Simulation results in Example~\ref{discont_initial_data} with fixed $\delta=0.2$ and two different values of $\beta$: 4.0 (first row) and 1.0 (second row) at four different time instances.}
    \label{heat_sim2}
\end{figure}

\begin{ex}\label{discont_initial_data}
	Let $\Omega$ be the domain contained within the curve given in \eqref{ring_shape_domain} and let
	\[ u_0(x,y)  = 
	\begin{cases}
		2, &\text{ if } x^2+y^2<0.3,\\
		\exp({-x^4-y^4}), &\text{ otherwise}.
	\end{cases}
	\]
	We solve the nonlocal diffusion equation 
	\begin{eqnarray}
		\begin{cases}
			u_t(x,y,t) = \Ldel u(x,y,t),& (x,y,t)\in \Omega\times (0,T)\\
			u(x,y,0) = u_0(x,y),&(x,y)\in \Omega\\
			u(x,y,t) = 0, & (x,y)\in \collar
		\end{cases}
	\end{eqnarray}  
	with $\delta=0.2$ and two values of $\beta$: $4.0$ and $1.0$ at four different time instances: $t=0.001$, $0.005$, $0.01$, and $0.03$ with time step size $\tau=2.5\cdot10^{-7}$.  We choose the 2D-FC parameters $d=4$, $M=5$, and $h=k_1=k_2=0.005$. 
\end{ex}
In Figure~\ref{heat_sim2}, we present the numerical solutions for a fixed $\delta=0.2$ and two different values of $\beta$: $4.0$ and $1.0$ at different final time. For the case $\beta=4.0$, which corresponds to the classical heat equation, we observe an instantaneous smoothing effect at $\partial\Omega$ compared to the nonlocal counterpart where $\beta=1.0$. 

\section{Discussion}\label{discussion}
 Fourier spectral methods are traditionally well-suited for periodic problems. In this work, we extend their applicability to nonlocal equations on bounded domains by constructing periodic extensions of the problem via Fourier continuation. This enables efficient computations while maintaining high-order convergence for smooth solutions and effectively addressing challenges associated with nonlocal operators near boundaries.
 
 We also conduct numerical experiments to investigate solution regularity for two-dimensional nonlocal Poisson and diffusion equations on bounded domains, with particular attention to discontinuities. Our simulations reveal that for integrable kernels, discontinuities may develop in the solution. In contrast, for singular kernels with $\beta>n+1$, the numerical solution appears to approximate a continuous function within the interior of the domain--an observation formalized in Conjecture~\ref{conjecture:continuity}.
 
 While our convergence analysis focuses on smooth solutions, we observe a decline in convergence order when discontinuities are present. Extending this approach to problems involving sharp interfaces or discontinuities, where limited regularity can reduce the effectiveness of spectral methods, remains an important direction for future research.


\clearpage
\bibliographystyle{acm}
\bibliography{refs}

@article{johansson2013mpmath,
  title={mpmath: a Python library for arbitrary-precision floating-point arithmetic (version 0.18)},
  author={Johansson, Fredrik and others},
  journal={Official Wiki Page available at: http://code.google.com/p/mpmath},
  year={2013}
}

@article{alali2021fourier,
  title={Fourier multipliers for nonlocal Laplace operators},
  author={Alali, Bacim and Albin, Nathan},
  journal={Applicable Analysis},
  volume={100},
  number={12},
  pages={2526--2546},
  year={2021},
  publisher={Taylor \& Francis}
}

@article{alali2020fourier,
  title={Fourier spectral methods for nonlocal models},
  author={Alali, Bacim and Albin, Nathan},
  journal={Journal of Peridynamics and Nonlocal Modeling},
  volume={2},
  pages={317--335},
  year={2020},
  publisher={Springer}
}

@article{bruno2022two,
  title={Two-dimensional Fourier continuation and applications},
  author={Bruno, Oscar P and Paul, Jagabandhu},
  journal={SIAM Journal on Scientific Computing},
  volume={44},
  number={2},
  pages={A964--A992},
  year={2022},
  publisher={SIAM}
}

@article{albin2014discrete,
  title={Discrete periodic extension using an approximate step function},
  author={Albin, Nathan and Pathmanathan, Sureka},
  journal={SIAM Journal on Scientific Computing},
  volume={36},
  number={2},
  pages={A668--A692},
  year={2014},
  publisher={SIAM}
}

@book{saad2003iterative,
  title={Iterative methods for sparse linear systems},
  author={Saad, Yousef},
  year={2003},
  publisher={SIAM}
}

@article{mustapha2023regularity,
  title={REGULARITY OF SOLUTIONS FOR NONLOCAL DIFFUSION EQUATIONS ON PERIODIC DISTRIBUTIONS},
  author={Mustapha, Ilyas and Alali, Bacim and Albin, Nathan},
  journal={Journal of Integral Equations and Applications},
  volume={35},
  number={1},
  pages={81--104},
  year={2023},
  publisher={Rocky Mountain Mathematics Consortium Tempe, AZ, USA}
}

@article{hu2012peridynamic,
  title={Peridynamic model for dynamic fracture in unidirectional fiber-reinforced composites},
  author={Hu, Wenke and Ha, Youn Doh and Bobaru, Florin},
  journal={Computer Methods in Applied Mechanics and Engineering},
  volume={217},
  pages={247--261},
  year={2012},
  publisher={Elsevier}
}

@article{bobaru2015cracks,
  title={Why do cracks branch? A peridynamic investigation of dynamic brittle fracture},
  author={Bobaru, Florin and Zhang, Guanfeng},
  journal={International Journal of Fracture},
  volume={196},
  number={1-2},
  pages={59--98},
  year={2015},
  publisher={Springer}
}

@article{gilboa2009nonlocal,
  title={Nonlocal operators with applications to image processing},
  author={Gilboa, Guy and Osher, Stanley},
  journal={Multiscale Modeling \& Simulation},
  volume={7},
  number={3},
  pages={1005--1028},
  year={2009},
  publisher={SIAM}
}

@article{lejeune2017modeling,
  title={Modeling tumor growth with peridynamics},
  author={Lejeune, Emma and Linder, Christian},
  journal={Biomechanics and modeling in mechanobiology},
  volume={16},
  number={4},
  pages={1141--1157},
  year={2017},
  publisher={Springer}
}

@article{isiet2021review,
  title={Review of peridynamic modelling of material failure and damage due to impact},
  author={Isiet, Mewael and Mi{\v{s}}kovi{\'c}, Ilija and Mi{\v{s}}kovi{\'c}, Sanja},
  journal={International Journal of Impact Engineering},
  volume={147},
  pages={103740},
  year={2021},
  publisher={Elsevier}
}

@article{dang2024regularity,
  title={Regularity of Solutions for the Nonlocal Wave Equation on Periodic Distributions},
  author={Dang, Thinh and Alali, Bacim and Albin, Nathan},
  journal={arXiv preprint arXiv:2408.00912},
  year={2024}
}

@article{jafarzadeh2020efficient,
  title={Efficient solutions for nonlocal diffusion problems via boundary-adapted spectral methods},
  author={Jafarzadeh, Siavash and Larios, Adam and Bobaru, Florin},
  journal={Journal of Peridynamics and Nonlocal Modeling},
  volume={2},
  pages={85--110},
  year={2020},
  publisher={Springer}
}

@article{du2016asymptotically,
  title={Asymptotically compatible Fourier spectral approximations of nonlocal Allen--Cahn equations},
  author={Du, Qiang and Yang, Jiang},
  journal={SIAM Journal on Numerical Analysis},
  volume={54},
  number={3},
  pages={1899--1919},
  year={2016},
  publisher={SIAM}
}

@article{du2017fast,
  title={Fast and accurate implementation of Fourier spectral approximations of nonlocal diffusion operators and its applications},
  author={Du, Qiang and Yang, Jiang},
  journal={Journal of Computational Physics},
  volume={332},
  pages={118--134},
  year={2017},
  publisher={Elsevier}
}

@article{slevinsky2018spectral,
  title={A spectral method for nonlocal diffusion operators on the sphere},
  author={Slevinsky, Richard Mika{\"e}l and Montanelli, Hadrien and Du, Qiang},
  journal={Journal of Computational Physics},
  volume={372},
  pages={893--911},
  year={2018},
  publisher={Elsevier}
}

@book{bobaru2016handbook,
  title={Handbook of peridynamic modeling},
  author={Bobaru, Florin and Foster, John T and Geubelle, Philippe H and Silling, Stewart A},
  year={2016},
  publisher={CRC press}
}

@article{de2017finite,
  title={Finite element implementation of a peridynamic pitting corrosion damage model},
  author={De Meo, Dennj and Oterkus, Erkan},
  journal={Ocean Engineering},
  volume={135},
  pages={76--83},
  year={2017},
  publisher={Elsevier}
}

@article{dias2017review,
  title={A review of crack propagation modeling using peridynamics},
  author={Dias, Jo{\~a}o Paulo and Bazani, M{\'a}rcio Antonio and Paschoalini, Amarildo Tabone and Barbanti, Luciano},
  journal={Probabilistic prognostics and health management of energy systems},
  pages={111--126},
  year={2017},
  publisher={Springer}
}

@article{diehl2015simulation,
  title={Simulation of wave propagation and impact damage in brittle materials using peridynamics},
  author={Diehl, Patrick and Schweitzer, Marc Alexander},
  journal={Recent trends in computational engineering-CE2014},
  pages={251--265},
  year={2015},
  publisher={Springer}
}

@inproceedings{foss2016differentiability,
  title={Differentiability and integrability properties for solutions to nonlocal equations},
  author={Foss, Mikil and Radu, Petronela},
  booktitle={New Trends in Differential Equations, Control Theory and Optimization: Proceedings of the 8th Congress of Romanian Mathematicians},
  pages={105--119},
  year={2016},
  organization={World Scientific}
}

@article{foss2018existence,
  title={Existence and regularity of minimizers for nonlocal energy functionals},
  author={Foss, Mikil D and Radu, Petronela and Wright, Cory},
  journal={Differential and Integral Equations},
  volume={31},
  number={11-12},
  pages={807–832},
  year={2018}
}

@article{hinds2012dirichlet,
  title={Dirichlet’s principle and wellposedness of solutions for a nonlocal p-Laplacian system},
  author={Hinds, Brittney and Radu, Petronela},
  journal={Applied Mathematics and Computation},
  volume={219},
  number={4},
  pages={1411--1419},
  year={2012},
  publisher={Elsevier}
}

@article{hu2017peridynamics,
  title={Peridynamics for fatigue life and residual strength prediction of composite laminates},
  author={Hu, YL and Madenci, Erdogan},
  journal={Composite Structures},
  volume={160},
  pages={169--184},
  year={2017},
  publisher={Elsevier}
}

@article{lipton2014dynamic,
  title={Dynamic brittle fracture as a small horizon limit of peridynamics},
  author={Lipton, Robert},
  journal={Journal of Elasticity},
  volume={117},
  pages={21--50},
  year={2014},
  publisher={Springer}
}

@article{mengesha2014nonlocal,
  title={Nonlocal constrained value problems for a linear peridynamic Navier equation},
  author={Mengesha, Tadele and Du, Qiang},
  journal={Journal of Elasticity},
  volume={116},
  number={1},
  pages={27--51},
  year={2014},
  publisher={Springer}
}

@article{mengesha2016characterization,
  title={Characterization of function spaces of vector fields and an application in nonlinear peridynamics},
  author={Mengesha, Tadele and Du, Qiang},
  journal={Nonlinear Analysis},
  volume={140},
  pages={82--111},
  year={2016},
  publisher={Elsevier}
}

@article{silling2000reformulation,
  title={Reformulation of elasticity theory for discontinuities and long-range forces},
  author={Silling, Stewart A},
  journal={Journal of the Mechanics and Physics of Solids},
  volume={48},
  number={1},
  pages={175--209},
  year={2000},
  publisher={Elsevier}
}

@article{silling2007peridynamic,
  title={Peridynamic states and constitutive modeling},
  author={Silling, Stewart A and Epton, M and Weckner, Olaf and Xu, Jifeng and Askari, E23481501120},
  journal={Journal of elasticity},
  volume={88},
  pages={151--184},
  year={2007},
  publisher={Springer}
}

@article{silling2017modeling,
  title={Modeling shockwaves and impact phenomena with Eulerian peridynamics},
  author={Silling, Stewart A and Parks, Michael L and Kamm, James R and Weckner, Olaf and Rassaian, Mostafa},
  journal={International Journal of Impact Engineering},
  volume={107},
  pages={47--57},
  year={2017},
  publisher={Elsevier}
}

@article{burkovska2024efficient,
  title={An efficient proximal-based approach for solving nonlocal Allen-Cahn equations},
  author={Burkovska, Olena and Mustapha, Ilyas},
  journal={arXiv preprint arXiv:2410.06455},
  year={2024}
}

@article{madenci2014coupling,
  title={Coupling of the peridynamic theory and finite element method},
  author={Madenci, Erdogan and Oterkus, Erkan and Madenci, Erdogan and Oterkus, Erkan},
  journal={Peridynamic theory and its applications},
  pages={191--202},
  year={2014},
  publisher={Springer}
}

@article{albin2011spectral,
  title={A spectral FC solver for the compressible Navier--Stokes equations in general domains I: Explicit time-stepping},
  author={Albin, Nathan and Bruno, Oscar P},
  journal={Journal of Computational Physics},
  volume={230},
  number={16},
  pages={6248--6270},
  year={2011},
  publisher={Elsevier}
}

@article{burkovska2020affine,
  title={Affine approximation of parametrized kernels and model order reduction for nonlocal and fractional Laplace models},
  author={Burkovska, Olena and Gunzburger, Max},
  journal={SIAM Journal on Numerical Analysis},
  volume={58},
  number={3},
  pages={1469--1494},
  year={2020},
  publisher={SIAM}
}

@article{press2007numerical,
  title={Numerical recipes 3rd edition},
  author={Press, William H and Teukolsky, Saul A and Vetterling, William T and Flannery, Brian P},
  journal={The art of scientific computing},
  volume={3},
  year={2007},
  publisher={Cambridge University Press Cambridge}
}

\end{document}